\documentclass[reqno,a4paper,10pt]{amsart}
\usepackage[a4paper,centering]{geometry}
\usepackage[utf8]{inputenc}
\usepackage[english]{babel}
\usepackage[T1]{fontenc}
\usepackage{newtxtext}
\usepackage{amssymb,amsmath,amsthm,amsfonts}
\usepackage{newtxmath}
\usepackage{mathrsfs,dsfont,bbm}
\usepackage[mathcal]{euscript}
\usepackage{mathtools}
\usepackage{stmaryrd}
\usepackage{xcolor}
\usepackage{microtype}
\usepackage{graphicx}
\usepackage{cite}
\usepackage[colorlinks,linkcolor=blue,citecolor=blue,urlcolor=black]{hyperref}
\usepackage{bookmark}

\newtheorem{thm}{Theorem}[section]

\newtheorem*{thm*}{Theorem}
\newtheorem{cor}{Corollary}[section]

\newtheorem{lem}{Lemma}[section]

\newtheorem{ass}{A}

\theoremstyle{definition}
\newtheorem{defn}{Definition}[section]

\theoremstyle{remark}
\newtheorem{rem}{Remark}[section]

\numberwithin{equation}{section}

\def\N{{\mathbb N}}

\def\R{{\mathbb R}}

\def\E{{\mathbb E}}

\def\Y{{\mathscr Y}}
\def\G{{\mathscr W}}

\def\cP{{\mathscr{P}}}
\def\norma #1{\left\lVert #1 \right\rVert}

\def\P{{\mathbb{P}}}

\def\de{{\rm d}}

\def\F{{\mathscr{F}}}

\def\M{{\mathcal{M}}}
\def\nm #1{\left\langle #1 \right\rangle}

\def\W{{\mathcal W}}
\def\BL{{\rm BL}}

\def\cS{{\mathfrak S}}
\def\cT{{\mathfrak V}}

\def\y{{\boldsymbol y}}

\def\vH{{\boldsymbol H}}

\def\by{{\boldsymbol Y}}
\def\bz{{\boldsymbol Z}}
\def\bmu{{\boldsymbol{\mu}}}

\def\w2{\stackrel{2}{\rightharpoonup}}

\def\strong2{\stackrel{2}{\rightarrow}}

\definecolor{viola}{rgb}{0.3,0,0.7}
\definecolor{ciclamino}{rgb}{0.5,0,0.5}
\definecolor{rosso}{rgb}{0.8,0,0}

\title[Backward--forward interacting systems]{Backward--forward interacting systems with probability-valued strategies and common noise}

\author[A. Melchor Hernandez]{Anderson Melchor Hernandez}
\address[A. Melchor Hernandez]{Visiting Researcher, The Fields Institute, 222 College St, Toronto, ON M5T 3J1}
\email{amelchor@fieldsinstitute.ca, amelchorh@unal.edu.co}

\date{\today}
\keywords{Backward stochastic differential equations, interacting agents,
measure-valued dynamics, Wasserstein distance.}
\subjclass[2020]{60H10, 49N80, 49J45, 93E20}

\begin{document}

\begin{abstract}
We study a system of interacting agents in which a backward stochastic equation is coupled to the forward evolution of a probability-valued strategy distribution. The agents interact through their empirical state--strategy law and are exposed to a common Brownian noise. Under Lipschitz, mass-conservation, and positivity assumptions, we establish well-posedness of the finite-agent system on a sufficiently short time interval. We then identify its large-population limit as a conditional mean-field backward--forward system and prove conditional propagation of chaos. When agent characteristics take finitely many values, we obtain order $N^{-1}$ bounds for the expected squared particle--copy error and the expected squared Wasserstein distance of order one between the empirical law and the conditional limiting law. An adaptive-strategy model with replicator dynamics illustrates the theory, accompanied by a numerical comparison of the two errors.
\end{abstract}

\maketitle
\tableofcontents

\section{Introduction}
Large systems of interacting agents in which each agent carries both an individual state and a probability distribution over labels have been studied in the context of multi-population dynamics and their
mean-field limits
\cite{MorandottiSolombrino2019,AmbrosioFornasierMorandottiSavare2021,
AlmiMorandottiSolombrino2021,AlmiDEramoMorandottiSolombrino2022}. The probability-valued component describes an agent's distribution over possible labels, rather than assigning it a fixed label; in the application considered here, these labels represent strategies.
Interactions may then influence both the individual state and the evolution of the label distribution. This raises two related questions: whether the finite-agent dynamics are well posed and whether they admit a large-population limit.

A forward stochastic formulation of the state--label setting was studied in \cite{DonofrioMelchorHernandez2026}. In that work, the
individual state follows a forward stochastic differential equation, while a probability measure describes the agent's distribution of labels. The analysis includes well-posedness and mean-field limits, with Eulerian and Lagrangian descriptions of the limiting dynamics. The present paper retains the state--probability structure but prescribes the individual state at the terminal time: the state satisfies a backward stochastic differential equation, whereas the probability-valued strategy component evolves forward from an initial distribution.

The difficulty is not merely the presence of a backward
equation alongside a forward probability-valued component.
The two unknowns have boundary data at opposite ends of the
time interval, and both coefficients depend on the empirical
state--strategy law. A fixed-point argument must therefore
control the backward state and its stochastic integrand
together with a forward evolution that remains in
$\mathcal P(U)$. In particular, freezing the input of the
backward equation does not suffice: the frozen
measure-valued equation must also be solved within
$\mathcal P(U)$, and the estimates for the two outputs must be combined with constants independent of the population size.

Evolving probability-valued strategies have been studied in
\cite{AmbrosioFornasierMorandottiSavare2021}. Mean-field BSDEs
obtained through a limit approach were studied in
\cite{BuckdahnDjehicheLiPeng2009}, while backward propagation
of chaos is addressed in
\cite{LauriereTangpi2022,PapapantoleonSaplaourasTheodorakopoulos,GaoHuWen}. In particular, \cite{PapapantoleonSaplaourasTheodorakopoulos} establishes well-posedness and propagation of chaos for general mean-field BSDEs, and \cite{GaoHuWen} obtains quantitative results for mean-field forward--backward systems. These results do not directly provide the coupled construction used here. Our forward
unknown is an agent-level process constrained to remain in $\mathcal P(U)$, and its evolution must be controlled together with the backward state and stochastic integrand. Likewise, the forward stochastic state--label system of \cite{DonofrioMelchorHernandez2026} has no backward state prescribed at the terminal time. The contribution here is a short-time well-posedness and conditional propagation-of-chaos argument for the coupled system \eqref{intro-system}. We construct its probability-valued component by a positivity- and mass-preserving iteration, derive population-uniform stability estimates for the coupled system, and compare the particles with copies conditioned on the common noise. When the marks have finite support, the conditional limiting paths reduce to finitely
many types, allowing the sampling error and the resulting
particle--copy error to be bounded at order $N^{-1}$ in
expected squared distance.

A backward formulation is natural when an agent's current state represents the value or assessment of a future outcome: the terminal quantity determines an earlier value \cite{ElKarouiPengQuenez1997}. In our model, the terminal outcome determines the backward state, while the strategy distribution is updated forward using information available at the current time. Moreover, all agents are driven by a common Brownian motion. Independently sampled agent characteristics do not eliminate this shared uncertainty as the population grows. As in other models with common noise, the relevant large-population description is therefore a conditional law, and independence of limiting agents is understood
conditionally on the common noise \cite{CarmonaZhu2016}. In the present setting, this law is the law of a representative agent conditional on the common Brownian motion.

More precisely, let $U$ be a compact metric space and set
$\Y=\R^d\times\cP(U)$, where $\cP(U)$ denotes the space of
probability measures on $U$. Agent $i$ has a backward state $Y^i$ and a
strategy distribution $\mu^i$. Writing
\begin{align}\label{empirical-strategy}
\nu_t^N=\frac1N\sum_{i=1}^N\delta_{(Y_t^i,\mu_t^i)}    
\end{align}
for the empirical state--strategy distribution, we study the
coupled system
\begin{align}\label{intro-system}
\begin{cases}
\displaystyle
Y_t^i=\xi^i+
\int_t^T H\bigl(Z_s^i,(Y_s^i,\mu_s^i),\nu_s^N\bigr)\,\de s
-\int_t^T Z_s^i\,\de B_s,\\[0.4em]
\displaystyle
\mu_t^i=\overline{\mu}_0^i+
\int_0^t F\bigl(Z_s^i,(Y_s^i,\mu_s^i),\nu_s^N\bigr)\,\de s,
\qquad i=1,\ldots,N.
\end{cases}
\end{align}
Here $B$ is the common Brownian motion and $Z^i$ is the
adapted process appearing as the integrand in the stochastic integral of the backward equation. It describes the martingale component of agent $i$'s backward state. The driver $H$ may depend on $Y^i$, $Z^i$, the agent's strategy distribution $\mu^i$, and the empirical law $\nu^N$. The vector field $F$ governs the forward evolution of the strategy distribution and may depend on the same quantities. Thus \eqref{intro-system} couples a terminal condition for $Y^i$ with an initial condition for $\mu^i$. 

The shared Brownian motion makes the population limit potentially random. The relevant question is therefore
whether $\nu^N$ converges to the law of a representative
state--strategy pair conditional on this common noise, and whether finite collections of agents become asymptotically independent after conditioning. This is the notion of conditional propagation of chaos used below
\cite{CarmonaZhu2016,Nikolaev2025}. 

Our first result establishes well-posedness on a sufficiently short time interval whose length is independent of the population size. It relies on Lipschitz conditions for $H$ and $F$, together with mass-conservation and positivity conditions that keep the forward component probability-valued.
\begin{thm*}[{Informal statement of \autoref{importantprop}}]
Under the Lipschitz, mass-conservation, and positivity assumptions
stated below, there exists $T_*>0$, independent of $N$, such that
for every $0<T\leq T_*$ and every $N\geq1$, system
\eqref{intro-system} has a unique solution.
\end{thm*}
The proof combines estimates for the frozen backward and
measure-valued equations in a contraction whose constants are independent of $N$.

To formulate the large-population limit, we assign each agent an independently sampled characteristic, or mark. A mark may specify the agent's terminal value $\xi^i$, its initial strategy distribution $\overline{\mu}_0^i$, or both. The marks are independent of the Brownian motion $B$ shared by all agents. They therefore provide agent-specific variation, whereas $B$ provides a source of uncertainty common to the population. Let $(\overline Y,\overline\mu)$ be the state--strategy
pair of a representative agent. At time $t$, its population law is taken conditional on the common-noise filtration
\begin{align*}
m_t=\mathcal L\bigl((\overline Y_t,\overline\mu_t)\mid\F_t^B\bigr).    
\end{align*}
For a fixed realization of the common noise up to time $t$, $m_t$ describes the distribution of state--strategy pairs across agents with different marks. As that realization changes, the distribution may change as well; thus $m_t$ is in general a random probability measure. In the limiting backward--forward system, $m_t$ replaces the empirical distribution $\nu_t^N$ in the interaction terms. Our second result identifies this conditional mean-field limit and establishes conditional propagation of chaos.

\begin{thm*}[{Informal statement of \autoref{thm:conditional-chaos}}]
Under suitable square-integrability assumptions on the terminal
data, there exists $T_{\mathrm{ch}}>0$ such that, for
$0<T\leq T_{\mathrm{ch}}$, the conditional mean-field system has
a unique solution. Moreover,
\begin{align*}
\E\sup_{t\in[0,T]}\W_1(\nu_t^N,m_t)^2
\longrightarrow 0,
\end{align*}
where $\W_{1}$ denotes the Wasserstein distance of order one. For each fixed number of agents, their limiting paths are independent conditionally on the common Brownian path.
\end{thm*}
The independence in the last statement concerns the limiting copies, not the agents in the finite system. At finite $N$, the agents interact through the empirical distribution $\nu^N$ and are therefore generally dependent even after conditioning on $B$. In the limiting construction, each copy has its own independently sampled mark, while all copies are driven by the same Brownian motion. Once the path of $B$ is fixed, the limiting copies are independent and have the same conditional law. For a fixed common-noise realization, averaging over the limiting copies approximates their conditional law $m_t$. The particle system can then be compared with those copies to control the difference between $\nu_t^N$ and $m_t$,
uniformly over the time interval. Without conditioning on $B$, the copies still share a source of randomness, so an ordinary law-of-large-numbers argument with a deterministic limiting law would not describe this limit. For populations with finitely many types, we also obtain a quantitative estimate. To state it, let $\mathcal{E}_N^2$ denote the average squared error between each particle and its corresponding limiting copy. Precisely, it is not a pointwise in time squared error but an aggregated path-space quantity. It combines the uniform in time squared differences between the backward states and between the strategy distributions with the time-integrated squared difference between the stochastic integrands, as made precise in \eqref{normaE}.
\begin{thm*}[{Informal statement of \autoref{thm:quantitative-chaos}}]
If, in addition, the distribution of the agent marks is
supported on at most $M$ points, then, for sufficiently small $T$, there are constants independent of $N$ such that
\begin{align*}
\mathcal E_N^2\leq \frac{C_1}{N},
\hskip 0.2cm \E\sup_{t\in[0,T]}\W_1(\nu_t^N,m_t)^2
\leq \frac{C_2}{N}.    
\end{align*}
The constants may depend on the coefficients, the time horizon, the number of types, and the second moment of the limiting agent path.
\end{thm*}
The empirical distribution of the limiting copies
therefore differs from $m_t$ only because the proportions of the finitely many types in the sample differ from their probabilities. These proportions have fluctuations of order $N^{-1/2}$, which yield an error of order $N^{-1}$ after squaring and taking expectations. The first estimate shows that the interacting particles remain close, on average, to their corresponding limiting copies at the same squared-error rate. Combining these two sources of error gives the second estimate for the empirical distribution, uniformly in time. The constants do not depend on the population size $N$. They may nevertheless depend on the number of types and on the properties of their limiting paths; the estimate is therefore a finite-type rate, rather than a rate uniform over arbitrary mark distributions. 

We illustrate the general theory with a population whose
agents adapt the weights assigned to finitely many strategies. Replicator dynamics provide a natural description of such adaptation: the weight of a strategy changes according to how its payoff compares with the agent's average payoff. This mechanism is familiar from evolutionary game dynamics
\cite{HofbauerSigmund2003} and has also been used in
mean-field models with evolving strategy distributions
\cite{AmbrosioFornasierMorandottiSavare2021}.

The remainder of the paper is organized as follows. Section \ref{sec:preliminaries} introduces the metric and functional framework and states the assumptions. Section \ref{sec:Results} presents the main results, whose proofs are given in Section \ref{sec:Proofs}. Section \ref{sec:application} develops an adaptive-strategy example and provides a numerical illustration of the results.

\section{Preliminaries}\label{sec:preliminaries}
\subsection{Probability measures}
Let $(X,d_X)$ be a metric space. We denote by $\M(X)$ the space of finite signed Borel measures on $X$, by $\M_+(X)$ its cone of
nonnegative measures, and by $\cP(X)$ the set of probability
measures. For $\sigma\in\M(X)$, its total variation measure is
denoted by $|\sigma|$. We also set
\begin{align*}
\M_0(X)\coloneqq \{\sigma\in\M(X):\sigma(X)=0\}. 
\end{align*}
The space of real-valued Lipschitz functions on $X$ is denoted by
${\rm Lip}(X,d_X)$. For $\phi\in{\rm Lip}(X,d_X)$, its Lipschitz
constant is
\begin{align*}
{\rm Lip}(\phi)\coloneqq \sup_{\substack{x,y\in X\\x\neq y}}
\frac{|\phi(x)-\phi(y)|}{d_X(x,y)}.    
\end{align*}
If $(X,d_X)$ is complete and separable, the $1$-Wasserstein
distance between $\mu_1,\mu_2\in\cP_1(X)$ is
\begin{align}\label{wass1}
\W_1(\mu_1,\mu_2)
\coloneqq
\inf_{\gamma\in\Gamma(\mu_1,\mu_2)}
\int_{X\times X}d_X(x_1,x_2)\,\gamma(\de x_1,\de x_2),
\end{align}
where $\Gamma(\mu_1,\mu_2)$ is the set of couplings of $\mu_1$ and $\mu_2$. Its Kantorovich--Rubinstein dual representation is

\begin{align*}
\W_1(\mu_1,\mu_2)=\sup_{\substack{\phi\in{\rm Lip}(X,d_X)\\{\rm Lip}(\phi)\leq1}}
\int_X\phi\,\de(\mu_1-\mu_2).    
\end{align*}

\begin{align}\label{spaceP1}
\cP_1(X) \coloneqq \left\{\mu\in\cP(X):
\int_X d_X(x,\overline x)\,\de\mu(x)<+\infty
\text{ for some }\overline x\in X\right\}.
\end{align}
This definition does not depend on the choice of $\overline{x}$. If $X$ is complete, then $(\cP_1(X),\W_1)$ is complete. Moreover, convergence in $\W_1$ is equivalent to weak convergence together with convergence of the first moments. Let $T>0$. We endow $C([0,T],X)$ with the uniform metric
\begin{align}\label{pathmetric}
d_{T,X}(\y_1,\y_2)\coloneqq 
\sup_{t\in[0,T]}d_X(\y_1(t),\y_2(t)).
\end{align}
\subsection{Our setting}\label{subsec2}
We study dynamics on $\Y\coloneqq \R^d\times\cP(U)$, where $(U,d_U)$ is a compact metric space. The component $x$ of $y=(x,\lambda)\in\Y$ is an agent's spatial state, while $\lambda$ is a probability distribution over the strategy space $U$. The space used to embed the strategy component is

\begin{align}\label{Arens-Ells}
\F(U)\coloneqq 
\overline{{\rm span}(\cP(U))}^{\,\norma{\cdot}_{\BL}}
\end{align}
where
\begin{align*}
\norma{\ell}_{\BL}
\coloneqq \sup_{\substack{\phi\in{\rm Lip}(U,d_U)\\
\norma{\phi}_{\rm Lip}\leq1}}\nm{\ell,\phi}    
\end{align*}
where we use the dual pairing $\nm{\ell,\phi}$ for $\ell\in\F(U)$ and $\phi\in{\rm Lip}(U,d_U)$,
and
\begin{align}\label{normLip}
\norma{\phi}_{\rm Lip}
\coloneqq \norma{\phi}_{\infty}+{\rm Lip}(\phi).
\end{align}
Since $U$ is compact, $\F(U)$ is a separable Banach space
containing $\M(U)$. For $\mu_1,\mu_2\in\cP(U)$, the two metrics satisfy

\begin{align}\label{disutile}
\norma{\mu_1-\mu_2}_{\BL}\leq\W_1(\mu_1,\mu_2)
\leq(1+D_U)\norma{\mu_1-\mu_2}_{\BL},
\end{align}

where
\begin{align}\label{costantDu}
D_U\coloneqq \min_{u_0\in U}\max_{u\in U}d_U(u_0,u)
\leq{\rm diam}(U).
\end{align}
In what follows, we set $\cS\coloneqq \R^d\times\F(U)$ and
\begin{align}
\norma{(x,\sigma)}_{\cS}\coloneqq |x|+\norma{\sigma}_{\BL}.    
\end{align}
We give $\Y\subset\cS$ the induced metric. Since $\cP(U)$
is compact in the $\BL$ metric, $\Y$ is a closed and
complete subset of $\cS$. For $\Psi\in\cP(\Y)$ and $1\leq p<+\infty$, set

\begin{align}\label{momentop}
m_p(\Psi)\coloneqq \displaystyle\int_\Y\norma{y}_{\cS}^{p}\,\de\Psi(y).    
\end{align}
In particular,
\begin{align}
\cP_1(\Y)=\{\Psi\in\cP(\Y):m_1(\Psi)<+\infty\}.    
\end{align}

Let $(\Omega,\F,(\F_t)_{t\in[0,T]},\P)$ be a complete
filtered probability space. In what follows, we use
\begin{align}
&L_T^2(\R^d)\coloneqq L^2(\Omega,\F_T,\P;\R^d),\\
&H_T^2(\R^{d\times d})\coloneqq \left\{Z\text{ predictable}:
\E\int_0^T|Z_t|^2\,\de t<+\infty
\right\},
\end{align}
where $|Z|$ is the Frobenius norm if confusion does not occur. In this work, for $i=1,\ldots,N$, we consider the system of forward-backward SDEs

\begin{align}\label{eq1}
\begin{cases}
\de Y_t^i=-H(Z_t^i,Y_t^i,\mu_t^i,\nu_t^N)\,\de t
+Z_t^i\,\de B_t,\\
\de\mu_t^i=F(Z_t^i,Y_t^i,\mu_t^i,\nu_t^N)\,\de t,\\
Y_T^i=\xi^i,\hskip 0.2cm \mu_0^i=\overline{\mu}_0^i.
\end{cases}
\end{align}
The maps
\begin{align}
H:\R^{d\times d}\times\Y\times\cP_1(\Y)\to\R^d,\hskip 0.2cm F:\R^{d\times d}\times\Y\times\cP_1(\Y)\to\F(U) 
\end{align}
are Borel measurable. Here

\begin{align}\label{empiricallaw}
\nu_t^N\coloneqq\frac1N\sum_{i=1}^N\delta_{(Y_t^i,\mu_t^i)}    
\end{align}
is the empirical measure associated to the system \eqref{eq1}. Here we denote by $B$ a $d$-dimensional Brownian motion $B$ and that $(\F_t)$ is the usual augmentation of its natural filtration.  Let us now give the notion of solution to the system \eqref{eq1}.
\begin{defn}\label{def1}
For a fixed $N\in\N$, a solution of \eqref{eq1} is a family
$(Y^i,Z^i,\mu^i)_{i=1}^N$ satisfying, for every $i$:
\begin{enumerate}
\item[(a)]
$Y^i$ is continuous and adapted, $Z^i\in H_T^2(\R^{d\times d})$, and $\mu^i$ is continuous, adapted, and $\cP(U)$-valued;
\item[(b)]
$\E\sup_{t\in[0,T]}|Y_t^i|^2<+\infty$ and
\begin{align*}
\E\int_0^T\left(|H(Z_t^i,Y_t^i,\mu_t^i,\nu_t^N)|^2+
\norma{F(Z_t^i,Y_t^i,\mu_t^i,\nu_t^N)}_{\BL}^2\right)\,\de t<+\infty;
\end{align*}
\item[(c)]
almost surely, simultaneously for all $t\in[0,T]$,
\begin{align*}
Y_t^i&=\xi^i+\int_t^T H(Z_s^i,Y_s^i,\mu_s^i,\nu_s^N)\,\de s
-\int_t^T Z_s^i\,\de B_s,\\
\mu_t^i&=\overline{\mu}_0^i
+\int_0^t F(Z_s^i,Y_s^i,\mu_s^i,\nu_s^N)\,\de s.
\end{align*}
\end{enumerate}
The second integral is a Bochner integral in $\F(U)$.
\end{defn}
Let us now state our assumptions.
\subsection{Assumptions}\label{sub:assumptions}
\begin{ass}\label{A1}
There is $L_H\geq0$ such that, for all $z^1,z^2\in\R^{d\times d}$, $y^1,y^2\in\Y$, and $\Psi^1,\Psi^2\in\cP_1(\Y)$,
\begin{align*}
|H(z^1,y^1,\Psi^1)-H(z^2,y^2,\Psi^2)|\leq L_H\left(|z^1-z^2|
+\norma{y^1-y^2}_{\cS}+\W_1(\Psi^1,\Psi^2)\right).
\end{align*}
Moreover, $|H(0,(0,\lambda_0),\delta_{(0,\lambda_0)})|$
is finite for every, $\lambda_0\in\cP(U)$.
\end{ass}
\begin{ass}\label{A2}
There is $L_F\geq0$ such that, for all
$z^1,z^2\in\R^{d\times d}$, $y^1,y^2\in\Y$, and
$\Psi^1,\Psi^2\in\cP_1(\Y)$,
\begin{align*}
\norma{F(z^1,y^1,\Psi^1)-F(z^2,y^2,\Psi^2)}_{\BL}
\leq L_F\left(|z^1-z^2|+\norma{y^1-y^2}_{\cS}+\W_1(\Psi^1,\Psi^2)
\right).
\end{align*}
Moreover,
$\norma{F(0,(0,\lambda_0),\delta_{(0,\lambda_0)})}_{\BL}$
is finite every $\lambda_0\in\cP(U)$.
\end{ass}
\begin{ass}\label{A3}
For all $(z,y,\Psi)\in\R^{d\times d}\times\Y\times\cP_1(\Y)$,
\begin{align*}
\nm{F(z,y,\Psi),1}=0.    
\end{align*}
\end{ass}
\begin{ass}\label{A4}
There exists $\delta\geq0$ such that, for all $z\in\R^{d\times d}$, $x\in\R^d$,$\lambda\in\cP(U)$, and $\Psi\in\cP_1(\Y)$,
the element
\begin{align*}
F(z,(x,\lambda),\Psi)+\delta\lambda   
\end{align*}
belongs to $\M_+(U)$.
\end{ass}
\begin{rem}
The conditions in \autoref{A1} and \autoref{A2} merely record that the coefficients take finite values. Together with Lipschitz continuity they imply the linear-growth estimates used below. On the other hand, \autoref{A3} expresses that perturbations of a probability measure occur through its tangent space. \autoref{A4} requires $F(z,(x,\lambda),\Psi)+\delta\lambda$ to be a
nonnegative measure; together with \autoref{A3},
this permits the probability-preserving construction later on.
\end{rem}

\begin{rem}\label{rem:integral-meaning}
Before proceeding, let us comment on the meaning of the integral
equation defining $\mu_t^i$ in \autoref{def1}(c). Since $F$ takes
values in the Banach space $\F(U)$, the identity
\begin{align*}
\mu_t^i=\overline{\mu}_0^i
+\int_0^t F(Z_s^i,Y_s^i,\mu_s^i,\nu_s^N)\,\de s
\end{align*}
is understood as an equality of $\F(U)$-valued random variables, where the integral is a Bochner integral. This requires only that
$s\mapsto F(Z_s^i,Y_s^i,\mu_s^i,\nu_s^N)$ be strongly measurable
and Bochner-integrable on $[0,T]$, which follows from
\autoref{A2} together with the growth estimate \eqref{boundF}
below: by \autoref{A2} the map is strongly measurable as a
composition of continuous and progressively measurable maps,
and by \eqref{boundF} and Jensen's inequality,
\begin{align*}
\E\left\|\int_s^t F(Z_r^i,Y_r^i,\mu_r^i,\nu_r^N)\,\de r
\right\|_{\BL}^2
\leq (t-s)\int_s^t
\E\,C_F^2\left(1+|Z_r^i|+\norma{(Y_r^i,\mu_r^i)}_{\cS}+m_1(\nu_r^N)\right)^2\,\de r,
\end{align*}
which is finite for all $0\leq s<t\leq T$ by the moment bounds established in Step 2 in the proof of \autoref{importantprop} below. In particular, the integral equation is understood in
$\F(U)$; the positivity of integrals appearing in the
construction below is justified by
\autoref{lem:positive-bochner}. The argument uses Bochner
integrability and does not assume that $\F(U)$ has the
Radon--Nikodym property; see also
\cite[Remark~2.2]{DonofrioMelchorHernandez2026}.
\end{rem}

\begin{lem}\label{lem:positive-bochner}
Let $U$ be a compact metric space. The cone $\M_+(U)$ is
closed in $\F(U)$ with respect to the $\BL$ norm. Consequently,
if $g:[0,T]\to\F(U)$ is Bochner integrable and
$g_s\in\M_+(U)$ for almost every $s$, then
\begin{align*}
\int_0^t g_s\,\de s\in\M_+(U)
\qquad\text{for every }t\in[0,T].
\end{align*}
\end{lem}
\begin{proof}
Let $(\sigma_n)_{n\geq1}\subset\M_+(U)$ converge to
$\ell\in\F(U)$ in the $\BL$ norm. Since the constant
function $1$ belongs to the class of test functions,
\begin{align*}
\sigma_n(U)=\langle\sigma_n,1\rangle
\longrightarrow\langle\ell,1\rangle.
\end{align*}
Thus $(\sigma_n(U))_n$ is bounded. By compactness of $U$,
a subsequence converges weakly to some $\sigma\in\M_+(U)$. For every Lipschitz function $\phi$, weak convergence and convergence in the $\BL$ norm give
\begin{align*}
\nm{\sigma,\phi}
=\lim_{n\to\infty}\nm{\sigma_n,\phi}
=\nm{\ell,\phi}.
\end{align*}
The embedding of finite measures into $\F(U)$ is injective, so $\ell=\sigma$. This proves that $\M_+(U)$ is closed in $\F(U)$.

For the integral assertion, use the construction of the
Bochner integral by approximation with simple functions,
see \cite[Appendix~A.2]{AmbrosioFornasierMorandottiSavare2021}.
Because $g$ is strongly measurable and takes values in the closed cone $\M_+(U)$ almost everywhere, it can be approximated in $L^1([0,t];\F(U))$ by simple functions taking values in $\M_+(U)$. Each of their integrals is a finite nonnegative linear combination of positive measures. These integrals converge in the $\BL$ norm to $\int_0^t g_s\,\de s$, which belongs to $\M_+(U)$ by the closedness just proved.
\end{proof}
The following Lemma is standard, but for the sake of completeness, we provide its proof \cite{PardouxPeng1990}.
\begin{lem}\label{lembase}
Let $h$ be a progressively measurable $\R^d$-valued process satisfying
\begin{align}
\E\int_0^T|h_t|^2\,\de t<+\infty,   
\end{align}
and let $\xi\in L_T^2(\R^d)$. Then
\begin{align}\label{eq2}
\begin{cases}
\de Y_t=-h_t\,\de t+Z_t\,\de B_t,\\
Y_T=\xi
\end{cases}
\end{align}
has a unique solution with $Y$ continuous and adapted,
$\E\sup_{t\in[0,T]}|Y_t|^2<+\infty$, and $Z\in H_T^2(\R^{d\times d})$. If $(Y^1,Z^1)$ and $(Y^2,Z^2)$ have the same terminal
value $\xi$ and drivers $h^1,h^2$, then
\begin{align*}
\E\sup_{t\in[0,T]}|Y_t^1-Y_t^2|^2
&\leq 4T\,\E\int_0^T|h_t^1-h_t^2|^2\,\de t,\\
\E\int_0^T|Z_t^1-Z_t^2|^2\,\de t
&\leq T\,\E\int_0^T|h_t^1-h_t^2|^2\,\de t.
\end{align*}
\end{lem}

\begin{proof}
Let us define
\begin{align}
\eta\coloneqq \xi+\int_0^T h_s\,\de s.    
\end{align}
The Cauchy--Schwarz inequality gives $\eta\in L_T^2(\R^d)$. Since $(\F_t)$ is the augmented Brownian filtration, the martingale representation theorem gives a unique $Z\in H_T^2(\R^{d\times d})$ such that
\begin{align}
M_t\coloneqq \E[\eta\mid\F_t]
=\E[\eta]+\int_0^t Z_s\,\de B_s.    
\end{align}
Let us define

\begin{align}
Y_t\coloneqq M_t-\int_0^t h_s\,\de s.    
\end{align}
Then $Y$ is continuous and adapted, $Y_T=\xi$, and the integrated form of \eqref{eq2} holds. Furthermore, by Doob's inequality and the Cauchy--Schwarz inequality yield that 

\begin{align}
\E\sup_{t\in[0,T]}|Y_t|^2<+\infty.
\end{align}
Now, let us consider the estimates, and write
$\Delta h=h^1-h^2$, $\Delta Y=Y^1-Y^2$, and
$\Delta Z=Z^1-Z^2$. Then
\begin{align*}
\Delta Y_t &= \E\left[\int_t^T\Delta h_s\,\de s
\,\middle|\,\F_t \right],\\
\Delta M_t&\coloneqq \Delta Y_t+\int_0^t\Delta h_s\,\de s
=\E\left[\int_0^T\Delta h_s\,\de s\,\middle|\,\F_t
\right].
\end{align*}
Hence, by conditional Jensen's inequality, Doob's $L^2$
inequality, and Cauchy--Schwarz one has
\begin{align*}
\E\sup_{t\in[0,T]}|\Delta Y_t|^2
&\leq\E\sup_{t\in[0,T]}\left(\E\left[\int_0^T|\Delta h_s|\,\de s
\,\middle|\,\F_t\right]\right)^2\\
&\leq4\E\left(\int_0^T|\Delta h_s|\,\de s\right)^2\\
&\leq 4T\,\E\int_0^T|\Delta h_s|^2\,\de s.
\end{align*}
To estimate $\Delta Z$, recall that the two solutions have the same terminal value. Set
\begin{align*}
A\coloneqq \int_0^T\Delta h_s\,\de s.
\end{align*}
For the difference of the martingales constructed above, we have
$\Delta M_t=\E[A\mid\F_t]$. In particular, $\Delta M_T=A$ and, since $\F_0$ is trivial up to null sets, $\Delta M_0=\E[A]$. The martingale representation therefore gives
\begin{align*}
\int_0^T\Delta Z_s\,\de B_s&=\Delta M_T-\Delta M_0=A-\E[A].
\end{align*}
Applying the It\"o isometry, then the variance identity and the
Cauchy--Schwarz inequality, yields
\begin{align*}
\E\int_0^T|\Delta Z_s|^2\,\de s
&=\E|A-\E[A]|^2=\E|A|^2-|\E[A]|^2\\
&\leq\E|A|^2=\E\left|\int_0^T\Delta h_s\,\de s\right|^2\\
&\leq T\,\E\int_0^T|\Delta h_s|^2\,\de s.
\end{align*}
Lastly, by using these estimates, we obtain uniqueness.
\end{proof}
In what follows, we consider
\begin{align}
\widehat{\G}^{N}\coloneqq (\R^{d\times d})^N\times(\R^d)^N
\times(\cP(U))^N, \hskip 0.2cm \widehat{\cT}^{N}\coloneqq 
(\R^{d\times d})^N\times(\R^d)^N
\times(\F(U))^N.   
\end{align}
We set $\bmu\coloneqq (\mu^1,\ldots,\mu^N)$, $\by\coloneqq (Y^1,\ldots,Y^N)$, and $\bz\coloneqq (Z^1,\ldots,Z^N)$. For
$\Psi\in\cP_1(\Y)$, set
\begin{align}\label{H_vectorialfields}
\begin{aligned}
&\vH^N(\bz,\by,\bmu,\Psi)\coloneqq 
\bigl(H(Z^1,Y^1,\mu^1,\Psi),\ldots,H(Z^N,Y^N,\mu^N,\Psi)\bigr),\\
&{\boldsymbol F}^N(\bz,\by,\bmu,\Psi)\coloneqq 
\bigl(F(Z^1,Y^1,\mu^1,\Psi),\ldots,F(Z^N,Y^N,\mu^N,\Psi)\bigr).
\end{aligned}
\end{align}
The stochastic integral with respect to the common Brownian
motion is defined componentwise:
\begin{align*}
\int_t^T\bz_s\,\de B_s\coloneqq \left(\int_t^T Z_s^1\,\de B_s,\ldots,\int_t^T Z_s^N\,\de B_s\right).
\end{align*}
Consequently \eqref{eq1} takes the form

\begin{align}\label{system2}
\begin{cases}
\displaystyle
\by_t=\boldsymbol\xi +\int_t^T\vH^N(\bz_s,\by_s,\bmu_s,\nu_s^N)\,\de s
-\int_t^T\bz_s\,\de B_s,\\
\displaystyle\bmu_t=\overline{\bmu}_0
+\int_0^t\boldsymbol F^N(\bz_s,\by_s,\bmu_s,\nu_s^N)\,\de s.
\end{cases}
\end{align}
Here ${\boldsymbol\xi}\coloneqq (\xi^1,\ldots,\xi^N)$ and
$\overline{\bmu}_0=(\overline{\mu}_0^1,\ldots,\overline{\mu}_0^N)$.
We assume ${\boldsymbol\xi}\in L^2(\Omega,\F_T,\P;\R^{dN})$
and $\overline{\bmu}_0\in(\cP(U))^N$.

\subsection{Agent marks and common-noise filtrations}
\label{subsec:marks}
 For the large-population analysis, we work on a probability space carrying a $d$-dimensional Brownian motion $B$ and an i.i.d. sequence $(\Theta^i)_{i\geq1}$ of random variables
taking values in a Polish space $E$. The marks are independent of $B$ and represent characteristics specific to each agent.
 
Let $G:E\times C([0,T],\R^d)\to\R^d$ and
$J:E\to\cP(U)$ be measurable maps. We specify the terminal
and initial data by
\begin{align}\label{chaos-data}
\xi^i=G(\Theta^i,B),\qquad
\overline\mu_0^i=J(\Theta^i),\qquad
\E|\xi^1|^2<+\infty.
\end{align}
Thus the same mark may determine both an agent's terminal value and its initial strategy distribution. The terminal value may also depend on the entire path of the common Brownian motion; it need not be measurable at time zero. For the particle system, we use the usual augmentation of
\begin{align}\label{chaos-filtration}
\F_t=\sigma\bigl(\Theta^i:i\geq1,\,
                 B_s:0\leq s\leq t\bigr).
\end{align}
We denote by $\F_t^B$ the usual augmentation of
$\sigma(B_s:0\leq s\leq t)$. For the limiting copy
associated with agent $i$, we use the smaller filtration
given by the usual augmentation of
\begin{align}\label{single-filtration}
\F_t^i=\sigma\bigl(\Theta^i,\,
                   B_s:0\leq s\leq t\bigr).
\end{align}
In particular, $\overline\mu_0^i$ is $\F_0^i$-measurable. These enlarged filtrations are used for the population-limit results. The well-posedness result below is first formulated for the
Brownian filtration introduced before Definition~\ref{def1}.
The enlarged filtrations defined here are used for the population-limit results.
\begin{lem}
\label{lem:independent-enlargement}
Let $\Theta$ be a random element taking values in a Polish space
and independent of the Brownian motion $B$. Let $(\mathcal G_t)_{t\in[0,T]}$
be the usual augmentation of
\begin{align*}
\mathcal G_t^0
\coloneqq \sigma(\Theta)\vee\sigma(B_s:0\leq s\leq t).
\end{align*}
Then $B$ is a Brownian motion with respect to $(\mathcal G_t)$.
Moreover, for every $\eta\in L^2(\Omega,\mathcal G_T,\P;\R^d)$,
there exists a unique $\mathcal G$-predictable process
$Z\in H_T^2(\R^{d\times d})$ such that
\begin{align*}
\E[\eta\mid\mathcal G_t]=\E[\eta\mid\mathcal G_0]
+\int_0^t Z_s\,\de B_s,\hskip 0.2cm 0\leq t\leq T.
\end{align*}
Consequently, \autoref{lembase} remains valid, with the same
two difference estimates, when its filtration is replaced
by $(\mathcal G_t)$.
\end{lem}
\begin{proof}
Independence of $\Theta$ and $B$ preserves the Brownian
property under this initial enlargement. The predictable
representation property of the Brownian filtration likewise
extends to an initial enlargement by independent information;
see \cite{Jacod1985,Protter2005}. This gives the displayed
representation, including its $\mathcal G_0$-measurable
initial value.

For completeness, we check the difference estimate that
changes notation when $\mathcal G_0$ is nontrivial. Let two
equations in \autoref{lembase} have the same terminal value
and drivers $h^1,h^2$, and put
\begin{align*}
A\coloneqq\int_0^T(h_s^1-h_s^2)\,\de s.
\end{align*}
The difference of their representing martingales equals
$\E[A\mid\mathcal G_t]$. Hence
\begin{align*}
\int_0^T(Z_s^1-Z_s^2)\,\de B_s
=A-\E[A\mid\mathcal G_0].
\end{align*}
It\"o's isometry, the $L^2$-contraction property of
conditional expectation, and Cauchy--Schwarz give
\begin{align*}
\E\int_0^T|Z_s^1-Z_s^2|^2\,\de s
=\E\bigl|A-\E[A\mid\mathcal G_0]\bigr|^2
\leq T\,\E\int_0^T|h_s^1-h_s^2|^2\,\de s.
\end{align*}
The estimate for $\E\sup_t|Y_t^1-Y_t^2|^2$ follows as in \autoref{lembase}.
\end{proof}
In what follows, we state our main results.
\section{Our results}\label{sec:Results}
Let us state our first result.
\begin{thm}\label{importantprop}
Assume that \autoref{A1}--\autoref{A4} hold true, and consider the usual augmentation of the Brownian filtration, with
$\boldsymbol\xi\in L^2(\Omega,\F_T,\P;\R^{dN})$ and
deterministic $\overline{\bmu}_0\in(\cP(U))^N$. There exists $T_{\ast}>0$, depending only on $L_H$ and $L_F$, such that, for every $N\in\N$ and every $0<T\leq T_\ast$, system \eqref{system2} has a unique solution in the sense
of \autoref{def1}.
\end{thm}
We now turn to the marked-agent setting of subsection~\ref{subsec:marks}. For the particle filtration
\eqref{chaos-filtration} and each single-agent filtration
\eqref{single-filtration}, the representation and difference
estimates needed below follow from \autoref{lem:independent-enlargement}. In particular, the time-zero conditional expectation need not be deterministic. In what follows, we denote by $(Y^{i,N},Z^{i,N},\mu^{i,N})_{i=1}^N$ the solution of the $N$-particle system with these data, and set
\begin{align}
X_t^{i,N}
&\coloneqq (Y_t^{i,N},\mu_t^{i,N}),\\
\nu_t^N
&\coloneqq
\frac1N\sum_{i=1}^N\delta_{X_t^{i,N}}.
\label{chaos-particle-empirical}
\end{align}

\begin{thm}[Conditional propagation of chaos]
\label{thm:conditional-chaos}
Assume \autoref{A1}--\autoref{A4} and \eqref{chaos-data}. There exists $T_{\mathrm{ch}}>0$, depending only on $L_H$ and $L_F$, such that, for every $0<T\leq T_{\mathrm{ch}}$, the following
statements hold. For every $i\geq1$, there exists a triple
$(\overline{Y}^i,\overline{Z}^i,\overline{\mu}^i)$ with the following properties:
\begin{enumerate}
\item[(a)] $\overline{Y}^i$ is a continuous, $\F^i$-adapted,
$\R^d$-valued process, and
\begin{align}
\E\sup_{t\in[0,T]}|\overline{Y}_t^i|^2<+\infty;
\end{align}

\item[(b)]
$\overline{Z}^i$ is an $\F^i$-predictable, $\R^{d\times d}$-valued process, and
\begin{align}
\E\int_0^T|\overline{Z}_t^i|^2\,\de t<+\infty;
\end{align}

\item[(c)]
$\overline{\mu}^i$ is a continuous, $\F^i$-adapted, $\cP(U)$-valued process.
\end{enumerate}
Set
\begin{align}
\label{chaos-conditional-law}
\overline{X}_t^i\coloneqq(\overline{Y}_t^i,\overline{\mu}_t^i),
\hskip 0.2cm m_t\coloneqq\mathcal L(\overline{X}_t^i\mid\F_t^B).
\end{align}
The conditional law $m$ does not depend on $i$. It admits an $\F^B$-adapted version with continuous paths in $(\cP_1(\Y),\W_1)$. For every $i$, the following equations hold almost surely, simultaneously for all $t\in[0,T]$:
\begin{align}\label{chaos-limit-backward}
\overline{Y}_t^i&=\xi^i+\int_t^TH\bigl(\overline Z_s^i,(\overline Y_s^i,\overline\mu_s^i),m_s\bigr)\,\de s-\int_t^T\overline Z_s^i\,\de B_s,
\\
\overline\mu_t^i&=\overline\mu_0^i+
\int_0^t F\bigl(\overline Z_s^i,(\overline Y_s^i,\overline\mu_s^i),m_s\bigr)\,\de s.
\label{chaos-limit-measure}
\end{align}
The drivers satisfy
\begin{align}
\E\int_0^T\Bigl(&\bigl|H\bigl(\overline Z_t^i,(\overline Y_t^i,\overline\mu_t^i),m_t\bigr)\bigr|^2\\
&+\norma{F\bigl(\overline Z_t^i,(\overline Y_t^i,\overline\mu_t^i),m_t\bigr)}_{\BL}^2\Bigr)\,\de t
<+\infty.
\end{align}
The solution is unique in this class: for each $i$, $\overline Y^i$ and $\overline\mu^i$ are unique up to indistinguishability, and $\overline Z^i$ is unique $\de t\otimes\de\P$-almost everywhere. The processes $(\overline X^i)_{i\geq1}$ are independent and identically distributed conditionally on the
entire Brownian path. For each $N$, couple particle $i$ with limit copy $i$ using the same pair $(\Theta^i,B)$. Then
\begin{align}\label{chaos-particle-convergence}
\begin{aligned}
\lim_{N\to\infty}\E\Bigg[\frac1N\sum_{i=1}^N
\Bigg(&\sup_{t\in[0,T]}|Y_t^{i,N}-\overline Y_t^i|^2\\
&+\int_0^T|Z_t^{i,N}-\overline Z_t^i|^2\,\de t\\
&+\sup_{t\in[0,T]}\norma{\mu_t^{i,N}-\overline\mu_t^i}_{\BL}^2
\Bigg)\Bigg]=0.
\end{aligned}
\end{align}
Moreover,
\begin{align}\label{chaos-empirical-convergence}
\lim_{N\to\infty}\E\sup_{t\in[0,T]}\W_1(\nu_t^N,m_t)^2=0.
\end{align}
In particular, for each fixed $k$, the first $k$ particle paths converge in probability to the first $k$ limit paths. The latter paths are independent conditionally on $B$.
\end{thm}

\begin{thm}[Quantitative conditional propagation of chaos]
\label{thm:quantitative-chaos}
Assume the hypotheses of \autoref{thm:conditional-chaos}, and suppose
in addition that the common law of the marks $\Theta^i$ is supported
on at most $M<+\infty$ points. Fix $\lambda_0\in\cP(U)$, let
$x_0\in C([0,T],\Y)$ be the constant path with value $(0,\lambda_0)$,
and set
\begin{align*}
K_T\coloneqq \E\,d_{T,\Y}(\overline X^1,x_0)^2,\hskip 0.2cm
C_{\mathrm{err}}&\coloneqq
40L_H^2+6L_F^2e^{2L_F}.
\end{align*}
Let
\begin{align*}
\overline\nu_t^N
\coloneqq \frac1N\sum_{i=1}^N\delta_{\overline X_t^i}, \hskip 0.2cm
\varepsilon_N\coloneqq
\E\sup_{t\in[0,T]}\W_1(\overline\nu_t^N,m_t)^2,
\end{align*}
and define
\begin{align}\label{normaE}
\begin{aligned}
\mathcal E_N^2\coloneqq\frac1N\sum_{i=1}^N\E\Bigg[
&\sup_{t\in[0,T]}|Y_t^{i,N}-\overline Y_t^i|^2
+\int_0^T|Z_t^{i,N}-\overline Z_t^i|^2\,\de t\\
&+\sup_{t\in[0,T]}
\norma{\mu_t^{i,N}-\overline\mu_t^i}_{\BL}^2
\Bigg].
\end{aligned}
\end{align}
Then $K_T<+\infty$ and, for every $N\geq1$ and $0<T\leq T_{\mathrm{ch}}$,
\begin{align}
\varepsilon_N &\leq \frac{M K_T}{N},
\label{chaos-quantitative-sampling}\\
\mathcal{E}_N^2 &\leq \frac{2C_{\mathrm{err}}T^2 M K_T}{N},
\label{chaos-quantitative-particles}\\
\E\sup_{t\in[0,T]}\W_1(\nu_t^N,m_t)^2&\leq \frac{(8C_{\mathrm{err}}T^2+2)M K_T}{N}.
\label{chaos-quantitative-empirical}
\end{align}
Moreover, for every fixed $k\leq N$,
\begin{align}\label{chaos-quantitative-fixed-k}
\begin{aligned}
\sum_{i=1}^k\E\Bigg[&\sup_{t\in[0,T]}|Y_t^{i,N}-\overline Y_t^i|^2+\int_0^T|Z_t^{i,N}-\overline Z_t^i|^2\,\de t\\
&+\sup_{t\in[0,T]}\norma{\mu_t^{i,N}-\overline\mu_t^i}_{\BL}^2
\Bigg]\leq \frac{2kC_{\mathrm{err}}T^2 M K_T}{N}.
\end{aligned}
\end{align}
\end{thm}

\section{proofs}\label{sec:Proofs}
\begin{proof}[{Proof of \autoref{importantprop}}]
Let us split the proof into several steps.

\smallskip
\noindent\textit{Step 1: Functional spaces}
\smallskip

Let us first specify the norms used for the $N$-agent system. Let us set $\boldsymbol a=(a^1,\ldots,a^N)$, and define 

\begin{align*}
|\boldsymbol a|_N^2\coloneqq \frac1N\sum_{i=1}^N|a^i|^2.    
\end{align*}
We use this notation both when $a^i\in\R^d$ and when
$a^i\in\R^{d\times d}$; in the latter case $|\cdot|$ is
the Frobenius norm. If $\boldsymbol\lambda=(\lambda^1,\ldots,\lambda^N)\in\F(U)^N$, we set

\begin{align*}
\norma{\boldsymbol\lambda}_{N,\BL}^2\coloneqq \frac1N\sum_{i=1}^N
\norma{\lambda^i}_{\BL}^2.    
\end{align*}
Let $\mathcal E_T$ be the set of triples $a=(\by,\bz,\bmu)$ with the following properties:
\begin{enumerate}
\item[(i)] $\by=(Y^1,\ldots,Y^N)$ is a continuous,
adapted, $(\R^d)^N$-valued process satisfying
\begin{align*}
\E\sup_{t\in[0,T]}|\by_t|_N^2<+\infty;
\end{align*}
\item[(ii)] $\bz=(Z^1,\ldots,Z^N)$ is a predictable,
$(\R^{d\times d})^N$-valued process satisfying
\begin{align*}
\E\int_0^T|\bz_t|_N^2\,\de t<+\infty;
\end{align*}
\item[(iii)] $\bmu=(\mu^1,\ldots,\mu^N)$ is a continuous,
adapted, $(\cP(U))^N$-valued process. Continuity here is
with respect to the $\BL$ metric $\norma{\cdot}_{N,\BL}$.
\end{enumerate}
Let us notice that no additional integrability condition is required in (iii). Indeed, every $\lambda\in\cP(U)$ satisfies
$\norma{\lambda}_{\BL}\leq1$ because for each $\norma{\phi}_{\rm Lip}\leq1$, one has 
\begin{align*}
|\nm{\lambda,\phi}|\leq\norma{\phi}_\infty\leq1.
\end{align*}
Consequently, for every $a\in\mathcal E_T$,
\begin{align*}
\sup_{t\in[0,T]}\norma{\bmu_t}_{N,\BL}^2\leq1
\qquad\text{almost surely}.
\end{align*}
In what follows, for $a=(\by,\bz,\bmu)$ and $\widetilde a=(\widetilde{\by},\widetilde{\bz},\widetilde{\bmu})$ in $\mathcal E_T$, define
\begin{align}
D_T(a,\widetilde a)^2
\coloneqq 
\E\Bigg[
&\sup_{t\in[0,T]}|\by_t-\widetilde{\by}_t|_N^2
+\int_0^T|\bz_t-\widetilde{\bz}_t|_N^2\,\de t\\
&+\sup_{t\in[0,T]}
\norma{\bmu_t-\widetilde{\bmu}_t}_{N,\BL}^2
\Bigg].
\end{align}
Let us notice that this quantity is finite, and its square root defines a metric on $\mathcal{E}_T$. Let us prove that $(\mathcal{E}_T,D_T)$ is complete. Let $(a^k)_{k\geq1}$ be a Cauchy sequence on $\mathcal{E}_{T}$, with $a^k=(\by^k,\bz^k,\bmu^k)$. The first components are Cauchy in the space of continuous adapted $(\R^d)^N$-valued processes with norm

\begin{align}
\left(\E\sup_{t\in[0,T]}|\by_t|_N^2\right)^{1/2}.    
\end{align}
There is therefore a continuous adapted process $\by$
such that
\begin{align}
\E\sup_{t\in[0,T]}|\by_t^k-\by_t|_N^2
\longrightarrow0.   
\end{align}
On the other hand, the second components are Cauchy in the
$L^2$ space of predictable $(\R^{d\times d})^N$-valued processes. Hence there is a predictable process $\bz$ for which

\begin{align}
\E\int_0^T|\bz_t^k-\bz_t|_N^2\,\de t\longrightarrow0.    
\end{align}
For the third component, regard each $\bmu^k$ initially
as a continuous adapted $\F(U)^N$-valued process. The space $\F(U)^N$ is Banach. Hence the continuous adapted $\F(U)^N$-valued processes with finite norm

\begin{align}
\left(
\E\sup_{t\in[0,T]}
\norma{\bmu_t}_{N,\BL}^2
\right)^{1/2}
\end{align}
form a Banach space. Consequently, there is a continuous
adapted $\F(U)^N$-valued process $\bmu$ such that

\begin{align}
\E\sup_{t\in[0,T]}\norma{\bmu_t^k-\bmu_t}_{N,\BL}^2
\longrightarrow0.    
\end{align}
It remains to check that this limit takes values in $(\cP(U))^N$. Since $U$ is compact, convergence in the $\BL$ norm of probability measures implies convergence in $\W_1$ by \eqref{disutile}. The space $(\cP(U),\W_1)$ is complete. Thus $\cP(U)$ is closed in $\F(U)$: a $\BL$-convergent sequence of probability measures has a limit that is again a probability measure. Applying this observation componentwise to the uniformly
convergent subsequence shows that, outside one null set,
$\mu_t^i\in\cP(U)$ for every $t\in[0,T]$ and every $i=1,\ldots,N$. Therefore $(\by,\bz,\bmu)\in\mathcal E_T$, and

\begin{align}
D_T(a^k,(\by,\bz,\bmu))\longrightarrow0.    
\end{align}
This proves that $(\mathcal E_T,D_T)$ is a complete metric space.

\smallskip
\noindent\textit{Step 2: Growth estimates.}
Fix $\lambda_0\in\cP(U)$, and set $y_0\coloneqq (0,\lambda_0)$, and $\Psi_0\coloneqq \delta_{y_0}$. Since $\norma{\lambda}_{\BL}\leq1,$ one gets $\norma{y_0}_{\cS}\leq1$. Let us now derive growth estimates for $H$ and $F$. For $\Psi\in\cP_1(\Y)$, the only coupling of $\Psi$
and $\Psi_0$ is $\Psi\otimes\delta_{y_0}$.
Consequently, by the triangle inequality,
\begin{align*}
\W_1(\Psi,\Psi_0)
&=
\int_\Y\norma{y-y_0}_{\cS}\,\de\Psi(y)\\
&\leq
\int_\Y
\left(\norma{y}_{\cS}+\norma{y_0}_{\cS}\right)
\,\de\Psi(y)\\
&\leq m_1(\Psi)+1.
\end{align*}
Let us set $h_0\coloneqq |H(0,y_0,\Psi_0)|$,and $f_0\coloneqq \norma{F(0,y_0,\Psi_0)}_{\BL}$. Using \autoref{A1} and the previous bound, we obtain
\begin{align}\label{boundH}
\begin{aligned}
|H(z,y,\Psi)|&\leq h_0+ L_H\left(|z|+\norma{y-y_0}_{\cS}+\W_1(\Psi,\Psi_0)\right)\\
&\leq h_0+L_H\left(|z|+\norma{y}_{\cS}+m_1(\Psi)+2\right)\\
&\leq C_H\left(1+|z|+\norma{y}_{\cS}+m_1(\Psi)\right),
\end{aligned}
\end{align}
where $C_H\coloneqq h_0+2L_H$. Similarly, \autoref{A2}
gives
\begin{align}\label{boundF}
\begin{aligned}
\norma{F(z,y,\Psi)}_{\BL}&\leq f_0+L_F\left(|z|+\norma{y-y_0}_{\cS}+\W_1(\Psi,\Psi_0)\right)\\
&\leq f_0+L_F\left(|z|+\norma{y}_{\cS}+m_1(\Psi)+2\right)\\
&\leq C_F\left(1+|z|+\norma{y}_{\cS}+m_1(\Psi)\right),
\end{aligned}
\end{align}
where $C_F\coloneqq f_0+2L_F$. Now fix $a=(\by,\bz,\bmu)\in\mathcal{E}_T$, and define
\begin{align}
y_t^i\coloneqq (Y_t^i,\mu_t^i), \hskip 0.2cm \nu_t^a\coloneqq 
\frac1N\sum_{i=1}^N\delta_{y_t^i}.    
\end{align}
For every $t$, the definition of $m_1$ and the bound 
$\norma{\mu_t^i}_{\BL}\leq1$ give

\begin{align}
\begin{aligned}
m_1(\nu_t^a)
&=
\frac1N\sum_{i=1}^N
\norma{(Y_t^i,\mu_t^i)}_{\cS}=
\frac1N\sum_{i=1}^N
\left(
|Y_t^i|+\norma{\mu_t^i}_{\BL}
\right)\\
&\leq
\frac1N\sum_{i=1}^N|Y_t^i|+1\leq
|\by_t|_N+1.
\end{aligned}
\end{align}
Then by plugging into \eqref{boundH}, and \eqref{boundF} one gets
\begin{align}\label{boundFH}
\begin{aligned}
|H(Z_t^i,y_t^i,\nu_t^a)|&\leq C_H\left(3+|Z_t^i|+|Y_t^i|+|\by_t|_N
\right),\\
\norma{F(Z_t^i,y_t^i,\nu_t^a)}_{\BL}&\leq C_F\left(
3+|Z_t^i|+|Y_t^i|+|\by_t|_N\right).
\end{aligned}
\end{align}
Let us now check the square-integrability explicitly for $F$ and $H$. Let us set For $G=H$ or $G=F$, let $C_G$ be the corresponding constant, and let $\norma{\cdot}_G$ mean $|\cdot|$ when $G=H$ and $\norma{\cdot}_{\BL}$ when $G=F$. Since $(b_1+b_2+b_3+b_4)^2\leq4\sum_{j=1}^4 b_j^2$, we have by \eqref{boundFH} that
\begin{align}
\begin{aligned}
\frac1N\sum_{i=1}^N \norma{G(Z_t^i,y_t^i,\nu_t^a)}_G^2
&\leq 4C_G^2\left(9+|\bz_t|_N^2+\frac1N\sum_{i=1}^N|Y_t^i|^2
+|\by_t|_N^2\right)\\
&=4C_G^2\left(9+|\bz_t|_N^2+2|\by_t|_N^2\right).
\end{aligned}
\end{align}
Therefore
\begin{align}\label{integrab_FH}
\begin{aligned}
&\E\int_0^T \frac1N\sum_{i=1}^N \norma{G(Z_t^i,y_t^i,\nu_t^a)}_G^2
\,\de t\\
&\qquad\leq 4C_G^2 \left(9T+\E\int_0^T|\bz_t|_N^2\,\de t
+2T\,\E\sup_{t\in[0,T]}|\by_t|_N^2\right)<+\infty.
\end{aligned}
\end{align}
Let us now take
$\widetilde a=(\widetilde{\by},\widetilde{\bz},
\widetilde{\bmu})\in\mathcal E_T$, and set $\widetilde y_t^i
=(\widetilde Y_t^i,\widetilde\mu_t^i)$. For every $t$, the measure
\begin{align*}
\gamma_t\coloneqq \frac1N\sum_{i=1}^N \delta_{(y_t^i,\widetilde y_t^i)}
\end{align*}

is a coupling of $\nu_t^a$ and
$\nu_t^{\widetilde a}$. Hence, by the definition
of $\W_1$,
\begin{align*}
\W_1(\nu_t^a,\nu_t^{\widetilde a})
&\leq
\int_{\Y\times\Y}
\norma{y-\widetilde y}_{\cS}\,
\de\gamma_t(y,\widetilde y)\\
&=
\frac1N\sum_{i=1}^N
\left(
|Y_t^i-\widetilde Y_t^i|
+
\norma{\mu_t^i-\widetilde\mu_t^i}_{\BL}
\right)\\
&\leq
|\by_t-\widetilde{\by}_t|_N
+
\norma{\bmu_t-\widetilde{\bmu}_t}_{N,\BL}.
\end{align*}
Let us set $R_t^2\coloneqq 
|\by_t-\widetilde{\by}_t|_N^2+\norma{\bmu_t-\widetilde{\bmu}_t}_{N,\BL}^2$, then $(a+b)^2\leq2(a^2+b^2)$ gives
\begin{align}\label{bound_empirics}
\W_1(\nu_t^a,\nu_t^{\widetilde a})^2
\leq2R_t^2.   
\end{align}
The preceding estimate will be used when comparing two inputs
of the Picard map because it controls the difference between their empirical measures by the distance between the corresponding processes. We also verify the time regularity of the empirical measure. Fix $a\in\mathcal E_T$ and $s,t\in[0,T]$. Consider the coupling
\begin{align}
\gamma_{s,t}^a\coloneqq \frac1N\sum_{i=1}^N\delta_{(y_s^i,y_t^i)}    
\end{align}
between $\nu_s^a$ and $\nu_t^a$. Therefore
\begin{align*}
\W_1(\nu_s^a,\nu_t^a)&\leq\frac1N\sum_{i=1}^N
\norma{y_s^i-y_t^i}_{\cS}\\
&=\frac1N\sum_{i=1}^N\left(|Y_s^i-Y_t^i|
+\norma{\mu_s^i-\mu_t^i}_{\BL}\right).
\end{align*}
Since outside a single null set, every $Y^i$ and $\mu^i$ has
a continuous path, then right-hand
side tends to zero as $t\to s$. Hence $t\mapsto\nu_t^a$ is almost surely continuous in $(\cP_1(\Y),\W_1)$.
With the previous bounds at hand, we define our Picard scheme.

\smallskip
\noindent\textit{Step 3: Definition of the Picard scheme.}
Fix an input $a=(\by,\bz,\bmu)\in\mathcal E_T$. We construct
$\Phi(a)=(\by',\bz',\bmu')$ component by component. Throughout this step, let $\nu^a$ be the empirical-measure process associated with the fixed input $a$.

\smallskip
\noindent\textit{Construction of $(\by',\bz')$.}
For each $i=1,\ldots,N$, set

\begin{align}
h_t^i\coloneqq H(Z_t^i,(Y_t^i,\mu_t^i),\nu_t^a).    
\end{align}
The measurability assumptions on $H$, together with the properties of the input processes and of $\nu^a$, imply that $h^i$ is progressively measurable. Moreover, \eqref{integrab_FH} gives

\begin{align}
\E\int_0^T|h_t^i|^2\,\de t<+\infty.
\end{align}
Since $\xi^i\in L^2(\Omega,\F_T,\P;\R^d)$, \autoref{lembase} yields a unique pair $(Y^{\prime i},Z^{\prime i})$ satisfying
\begin{align}\label{picard-backward}
Y_t^{\prime i}=\xi^i+\int_t^T h_s^i\,\de s-\int_t^T Z_s^{\prime i}\,\de B_s,\qquad t\in[0,T].
\end{align}
Here $Y^{\prime i}$ is continuous and adapted, and $Z^{\prime i}$ is predictable, with
\begin{align}
\E\sup_{t\in[0,T]}|Y_t^{\prime i}|^2<+\infty, \hskip 0.2cm
\E\int_0^T|Z_t^{\prime i}|^2\,\de t<+\infty.    
\end{align}
The equality in \eqref{picard-backward} holds almost surely simultaneously for all $t$.

\smallskip
\noindent\textit{Construction of $\bmu'$.}
For each $i$, we seek a continuous, adapted,
$\cP(U)$-valued solution of
\begin{align}\label{picard-measure}
\mu_t^{\prime i} =\overline{\mu}_0^i
+\int_0^t F(Z_s^i,(Y_s^i,\mu_s^{\prime i}),\nu_s^a)\,\de s.
\end{align}
In this equation $Z^i$, $Y^i$, and $\nu^a$ are fixed input processes. Fix $i$ for the remainder of the construction, and set

\begin{align}
z_s\coloneqq Z_s^i,\hskip 0.2cm x_s:=Y_s^i, \hskip 0.2cm \Psi_s\coloneqq \nu_s^a,\hskip 0.2cm \lambda_0:=\overline{\mu}_0^i.    
\end{align}
We must solve
\begin{align}\label{measure-equation-short}
\lambda_t=\lambda_0+\int_0^t F(z_s,(x_s,\lambda_s),\Psi_s)\,\de s.
\end{align}
To this aim let us take the constant $\delta\geq0$ from \autoref{A4} and define
\begin{align}\label{damped-iteration}
\begin{aligned}
\lambda_t^0&\coloneqq \lambda_0,\\
\lambda_t^{k+1}&\coloneqq e^{-\delta t}\lambda_0
+\int_0^t e^{-\delta(t-s)}\bigl[F(z_s,(x_s,\lambda_s^k),\Psi_s)
+\delta\lambda_s^k\bigr]\,\de s,\hskip 0.2cm k\geq0.
\end{aligned}
\end{align}
Let us check that the iteration is well-defined. Suppose that $\lambda^k$ is a continuous, adapted, $\cP(U)$-valued process. It is bounded in $\F(U)$ because $\norma{\lambda_s^k}_{\BL}\leq1$. Also, by \eqref{boundF}, and $\norma{(x_s,\lambda_s^k)}_{\cS}\leq|x_s|+1$ show that

\begin{align}
\E\int_0^T \norma{F(z_s,(x_s,\lambda_s^k),\Psi_s)}_{\BL}^2
\,\de s<+\infty.    
\end{align}
Thus the integrand in \eqref{damped-iteration} is time integrable almost surely. Its integral defines a continuous process. Since the integrand at time $s$ depends only on information available at time $s$, this process is adapted. The claim applies first to $\lambda^0$ and then inductively to every iterate. We next show, by induction, that each iterate remains
in $\cP(U)$. This is true for $\lambda^0=\lambda_0$. Suppose it is true for $\lambda^k$. By \autoref{A4}, the strongly measurable integrand
\begin{align*}
s\longmapsto
e^{-\delta(t-s)}
\bigl[F(z_s,(x_s,\lambda_s^k),\Psi_s)
      +\delta\lambda_s^k\bigr]
\end{align*}
takes values in $\M_+(U)$ for almost every $s\in[0,t]$.
It is Bochner integrable by the growth estimate established
above. Hence \autoref{lem:positive-bochner} shows that
its integral belongs to $\M_+(U)$. Since $e^{-\delta t}\lambda_0\in\M_+(U)$ as well, we conclude that $\lambda_t^{k+1}$ is a nonnegative measure. Moreover, \autoref{A3} gives $\nm{F(z_s,(x_s,\lambda_s^k),\Psi_s),1}=0$. Then, one gets
\begin{align}
\begin{aligned}
\nm{\lambda_t^{k+1},1}&=e^{-\delta t}\nm{\lambda_0,1}
+\int_0^t e^{-\delta(t-s)}\left[\nm{F(z_s,(x_s,\lambda_s^k),\Psi_s),1}+\delta\nm{\lambda_s^k,1}
\right]\de s\\
&=e^{-\delta t}+\delta\int_0^t e^{-\delta(t-s)}\,\de s
=1.
\end{aligned}
\end{align}
Hence $\lambda_t^{k+1}\in\cP(U)$ for every $t$, and the induction is complete. We now prove convergence of the iterates. For $k\geq1$, subtract the formulas for $\lambda^{k+1}$ and
$\lambda^k$. Since their terms involving $\lambda_0$
cancel,
\begin{align}
\begin{aligned}
\lambda_t^{k+1}-\lambda_t^k
=\int_0^t e^{-\delta(t-s)}\Bigl[&F(z_s,(x_s,\lambda_s^k),\Psi_s)
-F(z_s,(x_s,\lambda_s^{k-1}),\Psi_s)\\
&+\delta(\lambda_s^k-\lambda_s^{k-1})\Bigr]\,\de s.
\end{aligned}
\end{align}
By \autoref{A2}, one gets
\begin{align}\label{iteration-difference}
\norma{\lambda_t^{k+1}-\lambda_t^k}_{\BL}
\leq (L_F+\delta)\int_0^t \norma{\lambda_s^k-\lambda_s^{k-1}}_{\BL}\,\de s.
\end{align}
Both $\lambda_t^0$ and $\lambda_t^1$ are probability measures, so $\norma{\lambda_t^1-\lambda_t^0}_{\BL}\leq2$. Set $K\coloneqq L_F+\delta$. Starting from this bound and applying \eqref{iteration-difference} repeatedly gives, for every $k\geq0$ and $t\in[0,T]$,
\begin{align}\label{factorial-bound}
\norma{\lambda_t^{k+1}-\lambda_t^k}_{\BL}
\leq 2\frac{(Kt)^k}{k!}\leq 2\frac{(KT)^k}{k!}.
\end{align}
Because there are only countably many iterates, we can choose
one event of probability one on which all of them have continuous paths, take values in $\cP(U)$, and satisfy
\eqref{factorial-bound} for every $t$. Fix an outcome
in this event. For $m>n$, the triangle inequality gives

\begin{align}
\sup_{t\in[0,T]}\norma{\lambda_t^m-\lambda_t^n}_{\BL}
\leq\sum_{k=n}^{m-1}2\frac{(KT)^k}{k!}.    
\end{align}
The right-hand side tends to zero as $n\to\infty$. Thus $(\lambda^k)$ is uniformly Cauchy on $[0,T]$. Since $\F(U)$ is complete, there is a function $t\mapsto\lambda_t\in\F(U)$ such that
\begin{align}\label{uniform-measure-limit}
\sup_{t\in[0,T]}\norma{\lambda_t^k-\lambda_t}_{\BL}
\longrightarrow0.
\end{align}
We now verify the properties of $\lambda$. First, for the fixed outcome under consideration, $\lambda$ is the uniform limit of continuous $\F(U)$-valued functions; hence its path is
continuous. Second, for each fixed $t$, $\lambda_t$ is the pointwise limit of the $\F_t$-measurable variables $\lambda_t^k$, and is therefore $\F_t$-measurable. Third, for every $t$, $\lambda_t^k\in\cP(U)$ and $\lambda_t^k\to\lambda_t$ in $\F(U)$. Since $\cP(U)$ is closed in $\F(U)$, as established in Step 1, $\lambda_t\in\cP(U)$. These statements hold simultaneously for every $t$ on the chosen
event. On its null complement, define $\lambda_t\coloneqq \lambda_0$ for all $t$; with the completeness of the filtration, this gives a continuous, adapted, $\cP(U)$-valued process on the whole probability space. It remains to show that this process solves \eqref{measure-equation-short}. Define $g_s\coloneqq F(z_s,(x_s,\lambda_s),\Psi_s)$. The Step 2 growth estimate, applied to the probability-valued process $\lambda$, shows that
$g$ is time-integrable almost surely. Moreover,
\autoref{A2} and \eqref{uniform-measure-limit}
imply
\begin{align*}
&\sup_{t\in[0,T]}\norma{\int_0^t e^{-\delta(t-s)}
\Bigl[F(z_s,(x_s,\lambda_s^k),\Psi_s)+\delta\lambda_s^k-g_s-\delta\lambda_s\Bigr]\,\de s}_{\BL}\\
&\qquad\leq (L_F+\delta)T\sup_{s\in[0,T]}\norma{\lambda_s^k-\lambda_s}_{\BL}\longrightarrow0.
\end{align*}
Passing to the limit in \eqref{damped-iteration}
therefore gives
\begin{align}\label{limit-damped-equation}
\lambda_t=e^{-\delta t}\lambda_0 +\int_0^t e^{-\delta(t-s)}
(g_s+\delta\lambda_s)\,\de s,\hskip 0.2cm t\in[0,T].
\end{align}
Notice that \eqref{limit-damped-equation} is not yet the equation we set out to solve, it is its damped form. The integrand $g+\delta\lambda$ is time-integrable, so the right-hand side is an absolutely continuous $\F(U)$-valued function of
$t$. The variation-of-constants formula gives, for almost every $t$,
\begin{align}
\frac{\de}{\de t}\lambda_t =-\delta\lambda_t+(g_t+\delta\lambda_t)=g_t.    
\end{align}
Since \eqref{limit-damped-equation} also gives
$\lambda_{t=0}=\lambda_0$, integration yields
$$
\lambda_t
=
\lambda_0+\int_0^t
F(z_s,(x_s,\lambda_s),\Psi_s)\,\de s
\qquad\text{for every }t\in[0,T].
$$
This is exactly \eqref{measure-equation-short}. Finally, suppose that $\lambda$ and $\widetilde\lambda$ are two continuous, $\cP(U)$-valued solutions with the same initial measure and the same frozen processes $z,x,\Psi$. Subtracting their integral equations and applying \autoref{A2} gives

\begin{align}
\norma{\lambda_t-\widetilde\lambda_t}_{\BL}
\leq L_F\int_0^t\norma{\lambda_s-\widetilde\lambda_s}_{\BL}
\,\de s.   
\end{align}
Gronwall's inequality implies that the solutions coincide for every $t$, almost surely. We have therefore constructed the unique solution of \eqref{picard-measure}. Let us set $\mu^{\prime i}:=\lambda$. Repeating this construction for the finitely many indices $i=1,\ldots,N$ gives a continuous, adapted,
$(\cP(U))^N$-valued process $\bmu'$. The growth estimate for $F$ also gives, for each $i$,

\begin{align}
\E\int_0^T
\norma{
F(Z_t^i,(Y_t^i,\mu_t^{\prime i}),\nu_t^a)
}_{\BL}^2\,\de t<+\infty.   
\end{align}
Together with the properties of $(\by',\bz')$ proved above, this shows that $\Phi(a)\in\mathcal E_T$.

\smallskip
\noindent\textit{Step 4: Difference estimate for the backward
equations.}
Let $a,\widetilde{a}\in\mathcal{E}_T$ and use $\Delta$
for differences of corresponding processes. For each $i$,
write
\begin{align*}
\Delta h_t^i&\coloneqq H(Z_t^i,Y_t^i,\mu_t^i,\nu_t^a)
-H(\widetilde Z_t^i,\widetilde Y_t^i,\widetilde\mu_t^i,\nu_t^{\widetilde a}).
\end{align*}
By \autoref{A1},
\begin{align*}
|\Delta h_t^i|&\leq L_H\left(|\Delta Z_t^i|+|\Delta Y_t^i|
+\norma{\Delta\mu_t^i}_{\BL}+\W_1(\nu_t^a,\nu_t^{\widetilde a})
\right).
\end{align*}
The inequality $(b_1+b_2+b_3+b_4)^2\leq4\sum_{j=1}^4 b_j^2$, followed by the empirical-measure estimate \eqref{bound_empirics} in Step 2 , gives
\begin{align}
\begin{aligned}\label{bound_step4}
|\Delta\boldsymbol h_t|_N^2
&\leq 4L_H^2\left(|\Delta\bz_t|_N^2+|\Delta\by_t|_N^2
+\norma{\Delta\bmu_t}_{N,\BL}^2+\W_1(\nu_t^a,\nu_t^{\widetilde a})^2\right)\\
&\leq 4L_H^2\left(|\Delta\bz_t|_N^2+3R_t^2\right).
\end{aligned}
\end{align}
Since the two output BSDEs have the same terminal values, applying the two estimates of \autoref{lembase} componentwise and averaging over $i$ yields
\begin{align}
\begin{aligned}
&\E\left[\sup_{t\in[0,T]}|\Delta\by_t'|_N^2
+\int_0^T|\Delta\bz_t'|_N^2\,\de t\right]\\
&\qquad\leq 5T\,\E\int_0^T|\Delta\boldsymbol h_t|_N^2\,\de t\\
&\qquad\leq 20L_H^2T\left[\E\int_0^T|\Delta\bz_t|_N^2\,\de t
+3T\,\E\sup_{t\in[0,T]}R_t^2\right].
\end{aligned}
\end{align}
For every $T\leq1$, the last expression is at most $60L_H^2T\,D_T(a,\widetilde a)^2$.

\smallskip
\noindent\textit{Step 5: Difference estimate for the
measure-valued equations.}
For each $i$, subtract the two output equations and use
\autoref{A2}:
\begin{align}
\norma{\Delta\mu_t^{\prime i}}_{\BL}
&\leq L_F\int_0^t\left(|\Delta Z_s^i|
+|\Delta Y_s^i|+\norma{\Delta\mu_s^{\prime i}}_{\BL}
+\W_1(\nu_s^a,\nu_s^{\widetilde a})\right)\,\de s.
\end{align}
The integral form of Gronwall's inequality therefore gives
\begin{align}
\sup_{t\in[0,T]}\norma{\Delta\mu_t^{\prime i}}_{\BL}
&\leq L_F e^{L_FT}\int_0^T\left(|\Delta Z_s^i|
+|\Delta Y_s^i|+\W_1(\nu_s^a,\nu_s^{\widetilde a})
\right)\,\de s.
\end{align}
Cauchy--Schwarz in time and by $(b_1+b_2+b_3)^2\leq3(b_1^2+b_2^2+b_3^2)$ one gets
\begin{align}
\sup_{t\in[0,T]}\norma{\Delta\mu_t^{\prime i}}_{\BL}^2
&\leq 3L_F^2e^{2L_FT}T\int_0^T\left(|\Delta Z_s^i|^2
+|\Delta Y_s^i|^2+\W_1(\nu_s^a,\nu_s^{\widetilde a})^2
\right)\,\de s.
\end{align}
By taking expectation, and using \eqref{bound_empirics}, one has
\begin{align}
\begin{aligned}
&\E\sup_{t\in[0,T]}
\norma{\Delta\bmu_t'}_{N,\BL}^2\\
&\quad\leq 3L_F^2e^{2L_FT}T\left[\E\int_0^T|\Delta\bz_s|_N^2\,\de s+3T\,\E\sup_{t\in[0,T]}R_t^2\right].
\end{aligned}
\end{align}
For $T\leq1$, we have

\begin{align}
\E\sup_{t\in[0,T]}
\norma{\Delta\bmu_t'}_{N,\BL}^2\leq 9L_F^2e^{2L_F}T\,D_T(a,\widetilde a)^2.
\end{align}

\smallskip
\noindent\textit{Step 6: Contraction and identification of the
solution.}
Combining Steps 4 and 5 gives, whenever $T\leq1$,
\begin{align}
D_T\bigl(\Phi(a),\Phi(\widetilde a)\bigr)^2
&\leq \left(60L_H^2+9L_F^2e^{2L_F}\right)T\,D_T(a,\widetilde a)^2.
\end{align}
Choose $T_*>0$ such that $T_*\leq1$ and
\begin{align}
\left(60L_H^2+9L_F^2e^{2L_F}\right)T_*<1.
\end{align}
Then $\Phi$ is a strict contraction on $\mathcal{E}_T$ for each $0<T\leq T_*$. The Banach fixed-point theorem gives a unique
$a=(\by,\bz,\bmu)$ with $\Phi(a)=a$. At the fixed point, the frozen empirical measure is precisely $\nu_t^N$, the frozen backward coefficients use $(\by,\bz,\bmu)$, and the measure equations use the same triple. Thus the two defining equations of $\Phi$ become \eqref{system2}. The integrability established in Steps 2 and 3 gives condition (b) of \autoref{def1}. Finally,
every solution in the class of \autoref{def1} is a fixed
point of $\Phi$, by uniqueness of the two decoupled equations in Step 3. The fixed-point uniqueness therefore proves uniqueness of the system.
\end{proof}

\begin{cor}\label{cor:random-initial}
\autoref{importantprop} remains valid, with the same constant
$T_\ast$, if the Brownian filtration $(\F_t)$ is replaced by any
filtration obtained by initially enlarging it with an
independent random variable, and if $\overline{\bmu}_0$ is only
assumed to be $\F_0$-measurable and $(\cP(U))^N$-valued, rather
than deterministic.
\end{cor}
\begin{proof}
Apply \autoref{lem:independent-enlargement} with the
independent random element used to enlarge the filtration.
It provides the martingale representation and both difference estimates of \autoref{lembase}, with unchanged constants. The measure-valued iteration in Step 3 of the proof of \autoref{importantprop} is pathwise and requires only that $\overline{\bmu}_0$ be $\F_0$-measurable and $(\cP(U))^N$-valued. All other growth, completeness, and contraction estimates are unchanged. The same threshold $T_\ast$ therefore applies.
\end{proof}

In what follows we prove \autoref{thm:conditional-chaos}.
\begin{proof}[{Proof of \autoref{thm:conditional-chaos}}]
We divide the proof into four parts.

\smallskip
\noindent\textit{Step 1. Construction of the conditional limit.}
Consider triples $a=(Y,Z,\mu)$ adapted to the single-particle filtration $\F^1$, with $Y$ continuous, $Z$ predictable, $\mu$ continuous and $\cP(U)$-valued, and
\begin{align}
\E\sup_{t\in[0,T]}|Y_t|^2+\E\int_0^T|Z_t|^2\,\de t<+\infty.
\end{align}
Identify continuous processes up to indistinguishability
and predictable processes up to equality $\de t\otimes\de\P$-almost everywhere. Let us equip this space with
\begin{align}\label{chaos-single-metric}
\begin{aligned}
\mathfrak D_T(a,\widetilde a)^2\coloneqq
\E\Bigg[&\sup_{t\in[0,T]}|Y_t-\widetilde Y_t|^2
+\int_0^T|Z_t-\widetilde Z_t|^2\,\de t\\
&+\sup_{t\in[0,T]}\norma{\mu_t-\widetilde\mu_t}_{\BL}^2
\Bigg].
\end{aligned}
\end{align}
This is a complete metric space by the same argument as in Step 1 of the proof of \autoref{importantprop}. Given such an input $a$, let
\begin{align}\label{chaos-input-law}
m_t^a \coloneqq \mathcal L\bigl((Y_t,\mu_t)\mid\F_t^B\bigr).
\end{align}
We make precise the relation between conditioning on the entire
Brownian path and conditioning on the information available at
time $t$. Set $\mathsf X\coloneqq C([0,T],\Y)$, equipped with
the uniform metric $d_{T,\Y}$ from \eqref{pathmetric}, and
regard $X^a\coloneqq(Y,\mu)$ as an $\mathsf X$-valued random
variable. Let
\begin{align}
Q_B^a\coloneqq\mathcal L(X^a\mid\mathcal F_T^B)
\end{align}
be a regular conditional law. For $t\in[0,T]$, define the
evaluation map
\begin{align}
\operatorname{ev}_t:\mathsf X\to\Y,
\qquad \operatorname{ev}_t(x)\coloneqq x(t).
\end{align}
The corresponding image measure is defined, for every Borel
set $A\subseteq\Y$, by
\begin{align}
\widetilde m_t^a(A)
\coloneqq
Q_B^a\bigl(\{x\in\mathsf X:x(t)\in A\}\bigr).
\end{align}
We claim that $\widetilde m_t^a$ is a version of the
conditional law $m_t^a$ defined in \eqref{chaos-input-law}.

Indeed, $X_t^a$ is measurable with respect to $\sigma(\Theta^1)\vee\mathcal F_t^B$. The independence of
$\Theta^1$ and $B$, together with the independence of Brownian increments, implies that increments after $t$ are independent of this sigma-field. Thus, for every bounded Borel function $\varphi:\Y\to\R$,
\begin{align}
\int_\Y\varphi(y)\,\widetilde m_t^a(\de y)
&=\E[\varphi(X_t^a)\mid\mathcal F_T^B]\notag\\
&=\E[\varphi(X_t^a)\mid\mathcal F_t^B]
\qquad\text{almost surely}.
\end{align}
Since $\Y$ is Polish, this identifies $\widetilde m_t^a$
with $m_t^a$ for each fixed $t$. In particular, $\widetilde m_t^a$ is $\mathcal F_t^B$-measurable. We use this version henceforth and denote it by $m_t^a$. This construction also gives path continuity. For any $s,t\in[0,T]$, the law of $(x(s),x(t))$ under $Q_B^a$
is a coupling of $m_s^a$ and $m_t^a$. Therefore,
\begin{align}\label{chaos-law-time}
\W_1(m_s^a,m_t^a)\leq \int_{\mathsf X}d_\Y(x(s),x(t))\,Q_B^a(\de x).
\end{align}
Fix a constant path $x_0\in\mathsf X$. The assumptions
on $Y$ and the fact that $\mu$ is probability-valued give
\begin{align}
\int_{\mathsf X}d_{T,\Y}(x,x_0)\,Q_B^a(\de x)
&=\E[d_{T,\Y}(X^a,x_0)\mid\mathcal F_T^B]
<\infty
\qquad\text{almost surely}.
\end{align}
Outside this single null set, each path $x\in\mathsf X$
is continuous, and
\begin{align}
d_\Y(x(s),x(t))\leq 2d_{T,\Y}(x,x_0).
\end{align}
Dominated convergence under $Q_B^a$ shows that the right-hand side of \eqref{chaos-law-time} tends to zero whenever $s\to t$. Thus $m^a$ has continuous paths in $(\cP_1(\Y),\W_1)$. Since each time-$t$ marginal is $\mathcal F_t^B$-measurable, this continuous version is $\F^B$-adapted. Let us define a map $\mathcal T$ on these inputs as follows. Let use the frozen driver
\begin{align}
h_t^a\coloneqq H\bigl(Z_t,(Y_t,\mu_t),m_t^a\bigr)
\end{align}
to solve the backward equation with terminal value $\xi^1$. Next, with the same frozen $Z$, $Y$, and $m^a$, solve the probability-valued equation 
\begin{align}
\mu_t'=\overline\mu_0^1+\int_0^t
F\bigl(Z_s,(Y_s,\mu_s'),m_s^a\bigr)\,\de s.
\end{align}
The growth estimates of Step 2 and the constructions
of Step 3  in the proof of \autoref{importantprop} apply to these equations. They show that $\mathcal T(a)=(Y',Z',\mu')$ belongs to the same
complete metric space. Notice that $\overline\mu_0^1=J(\Theta^1)$ is $\F_0^1$-measurable. Furthermore, the positivity and mass calculations of Step 3 in the proof of \autoref{importantprop} hold pathwise with this initial measure. Let us instead verify the contraction estimate. Let $a=(Y,Z,\mu)$ and $\widetilde a=(\widetilde Y,\widetilde Z,
\widetilde\mu)$, and set
\begin{align}
r_t^2 \coloneqq |Y_t-\widetilde Y_t|^2+
\norma{\mu_t-\widetilde\mu_t}_{\BL}^2.
\end{align}
Conditionally on $\F_t^B$, the joint law of $(Y_t,\mu_t)$ and
$(\widetilde Y_t,\widetilde\mu_t)$ is a coupling of $m_t^a$ and $m_t^{\widetilde a}$. Thus,
\begin{align}
\label{chaos-conditional-coupling}
\begin{aligned}
\W_1(m_t^a,m_t^{\widetilde a})
&\leq \E\bigl[|Y_t-\widetilde Y_t|+\norma{\mu_t-\widetilde\mu_t}_{\BL}\mid\F_t^B\bigr],\\
\E\W_1(m_t^a,m_t^{\widetilde a})^2&\leq2\E r_t^2.
\end{aligned}
\end{align}
Substituting \eqref{chaos-conditional-coupling} into the
backward and measure estimates of Steps 4 and 5
of the proof of \autoref{importantprop} gives,
for $T\leq1$,
\begin{align}\label{chaos-limit-contraction}
\mathfrak D_T\bigl(\mathcal T(a),\mathcal T(\widetilde a)
\bigr)^2
\leq C_{\mathrm{fix}}T\,
\mathfrak D_T(a,\widetilde a)^2,
\end{align}
where
\begin{align}\label{costCfix}
C_{\mathrm{fix}} \coloneqq 60L_H^2+9L_F^2e^{2L_F}.
\end{align}
For sufficiently small $T$, Banach's theorem gives a unique fixed point $(\overline Y^1,\overline Z^1,\overline\mu^1)$. 
Writing $m_t=m_t^{\overline a}$ for this fixed point
$\overline{a}=(\overline Y^1,\overline Z^1,\overline\mu^1)$, the fixed-point equations are \eqref{chaos-limit-backward} and \eqref{chaos-limit-measure} for $i=1$. The limiting state--strategy path is a measurable function of the agent's mark and the common Brownian path. Indeed, start the fixed-point iteration from a deterministic triple. At each iteration, the input conditional law is determined by the Brownian path, while the backward and measure-valued steps use
only that path and the agent's own mark. Thus every
state--strategy iterate is measurable with respect to
$\sigma(\Theta^1,B)$. The iterates converge in the metric \eqref{chaos-single-metric}; in particular, their
state--strategy paths converge in probability in
$C([0,T],\Y)$. A subsequence converges almost surely, so the limit path has a measurable representation
$\overline X^1=\Phi(\Theta^1,B)$ for some measurable map
$\Phi:E\times C([0,T],\R^d)\to C([0,T],\Y)$. 

For every $i$, the same construction with $\Theta^i$ in
place of $\Theta^1$ yields $\overline X^i=\Phi(\Theta^i,B)$ almost surely. Uniqueness
shows that these copies use the same conditional measure
process $m$. Since the marks are i.i.d. and independent of $B$, the paths $(\overline X^i)_{i\geq1}$ are i.i.d.
conditionally on $B$, with time-$t$ conditional marginal
$m_t$.

\smallskip
\noindent\textit{Step 2. Empirical measure of the limit copies.}
Set
\begin{align}\label{chaos-copy-empirical}
\overline\nu_t^N\coloneqq
\frac1N\sum_{i=1}^N\delta_{\overline X_t^i},
\hskip 0.2cm
e_N(t)\coloneqq \W_1(\overline\nu_t^N,m_t).
\end{align}
We show that
\begin{align}\label{chaos-sampling-error}
\varepsilon_N \coloneqq\E\sup_{t\in[0,T]}e_N(t)^2
\longrightarrow0.
\end{align}
We begin with almost-sure convergence of an empirical
measure on path space. Recall that $\mathsf{X}=C([0,T],\Y)$ and is equipped with the uniform metric $d_{T,\Y}$ from \eqref{pathmetric}. Since $\Y$ is complete and separable, $\mathsf X$ is a Polish space. Fix $\lambda_0\in\cP(U)$, and denote by $x_0\in\mathsf X$ the constant path with value $(0,\lambda_0)$. Let us regard each
\begin{align}
\overline X^i \coloneqq (\overline Y_t^i,\overline\mu_t^i)_{t\in[0,T]}
\end{align}
as an $\mathsf{X}$-valued random variable. The continuity of its paths makes this well-defined. Moreover,
\begin{align}
\begin{aligned}
d_{T,\Y}(\overline X^1,x_0)&=\sup_{t\in[0,T]}
\left(|\overline Y_t^1|+\norma{\overline\mu_t^1-\lambda_0}_{\BL}
\right)\\
&\leq\sup_{t\in[0,T]}|\overline Y_t^1|+2,
\end{aligned}
\end{align}
because both $\overline\mu_t^1$ and $\lambda_0$ are
probability measures. It follows from the square-integrability
assumption in the theorem that
\begin{align}\label{chaos-path-moment}
\E\,d_{T,\Y}(\overline X^1,x_0)^2&\leq 2\E\sup_{t\in[0,T]}|\overline Y_t^1|^2+8<+\infty.
\end{align}
Let
\begin{align}
\mathcal{B}\coloneqq \sigma(B_s:0\leq s\leq T).
\end{align}
Since $\mathsf{X}$ is Polish, the conditional law of $\overline{X}^1$ given $\mathcal{B}$ admits a regular
version. Denote it by
\begin{align}
Q_B\coloneqq \mathcal L(\overline X^1\mid\mathcal B),
\end{align}
and define the empirical measure of the limit
paths by
\begin{align}
Q_B^N\coloneqq \frac1N\sum_{i=1}^N\delta_{\overline X^i}.
\end{align}
The notation $Q_B^N$ does not mean that this
empirical measure is measurable with respect to
$B$ alone because it also depends on the individual marks
$\Theta^1,\ldots,\Theta^N$. Since the marks are i.i.d. and independent of $B$, for almost every realization of $B$, their conditional joint law is the product law
\begin{align}
\label{chaos-conditional-product}
\mathcal L\bigl(\overline X^1,\overline X^2,\ldots\mid\mathcal B
\bigr)=Q_B^{\otimes\N}.
\end{align}
In particular, conditionally on $B$, the variables
$\overline X^i$ are i.i.d. with common law $Q_B$. Let us now check the first-moment condition on $Q_B$. By the defining property of conditional expectation,
\begin{align}
\int_{\mathsf{X}}d_{T,\Y}(x,x_0)\,Q_B(\de x)
=\E\bigl[d_{T,\Y}(\overline X^1,x_0)\mid\mathcal B
\bigr]<+\infty
\end{align}
almost surely. Thus $Q_B\in\cP_1(\mathsf{X})$ almost surely. Also, each $Q_B^N$ is a finite empirical measure and therefore belongs to $\cP_1(\mathsf{X})$. For completeness, we give the two steps that imply convergence in $\W_1^{\mathsf X}$. Because
$\mathsf X$ is Polish, there is a countable family $(f_j)_{j\geq1}$ of bounded continuous functions that determines weak convergence of probability measures on $\mathsf X$. Fix a Brownian path for which \eqref{chaos-conditional-product} holds and $Q_B$ has finite first moment. Applying the strong law of large numbers under the corresponding conditional probability, simultaneously to this countable family, gives
\begin{align}
\int_{\mathsf X}f_j(x)\,Q_B^N(\de x)
=\frac1N\sum_{i=1}^N f_j(\overline X^i)\longrightarrow
\int_{\mathsf X}f_j(x)\,Q_B(\de x)
\end{align}
for every $j$, with conditional probability one.
Hence
\begin{align}
Q_B^N\rightharpoonup Q_B
\end{align}
weakly, with conditional probability one. Apply the strong law once more to the integrable function $x\mapsto d_{T,\Y}(x,x_0)$. It gives
\begin{align}
\int_{\mathsf X}d_{T,\Y}(x,x_0)\,Q_B^N(\de x)
&=\frac1N\sum_{i=1}^Nd_{T,\Y}(\overline X^i,x_0)\\
&\longrightarrow\int_{\mathsf X}d_{T,\Y}(x,x_0)\,Q_B(\de x)
\end{align}
with conditional probability one. Weak convergence together with convergence of these first moments is equivalent to convergence in $\W_1$ on the Polish space $\mathsf X$. Therefore, for almost
every Brownian path,
\begin{align}
\W_1^{\mathsf X}(Q_B^N,Q_B)
\longrightarrow0
\end{align}
with conditional probability one. Integrating this conditional probability-one statement over the law of $B$, we conclude that
\begin{align}\label{chaos-path-lln}
\W_1^{\mathsf X}(Q_B^N,Q_B)
\longrightarrow0
\qquad\text{almost surely}.
\end{align}
Here $\W_1^{\mathsf X}$ denotes the Wasserstein distance associated with the uniform path metric $d_{T,\Y}$. On the other hand, since the time-$t$ marginal of $Q_B$ is $m_t$, it follows that
\begin{align}\label{chaos-evaluation-bound}
\sup_{t\in[0,T]}e_N(t)\leq\W_1^{\mathsf X}(Q_B^N,Q_B).
\end{align}
For clarity, let us also justify passage to expectations
of squares. Let
\begin{align}
V_i\coloneqq d_{T,\Y}(\overline X^i,(0,\lambda_0)),
\hskip 0.2cm q_B \coloneqq \E[V_1\mid B].
\end{align}
The triangle inequality through the reference path
gives
\begin{align}\label{chaos-path-domination}
\W_1^{\mathsf X}(Q_B^N,Q_B)\leq \frac1N\sum_{i=1}^N V_i+q_B.
\end{align}

The squares of the quantities on the right-hand side of
\eqref{chaos-path-domination} are uniformly integrable.
Indeed, set $W_i\coloneqq V_i^2$. By \eqref{chaos-path-moment}, $W_1$ is integrable. The de la Vall\'ee Poussin criterion provides a convex increasing function $\Phi:[0,\infty)\to[0,\infty)$ such that $\Phi(r)/r\to\infty$ as $r\to\infty$ and $\E\Phi(W_1)<\infty$; see
\cite[Theorem~3.36]{advanced-probability}. By Jensen's
inequality and the identical distribution of the $W_i$,
\begin{align}
\E\Phi\left(\frac1N\sum_{i=1}^N W_i\right)
\leq \frac1N\sum_{i=1}^N\E\Phi(W_i)
=\E\Phi(W_1).
\end{align}
Conditional Jensen's inequality also gives
\begin{align}
\E\Phi\bigl(\E[W_1\mid B]\bigr)\leq \E\Phi(W_1).
\end{align}
Thus the family
$\bigl\{N^{-1}\sum_{i=1}^N V_i^2\bigr\}_{N\geq1}$
is uniformly integrable. Since
$q_B=\E[V_1\mid B]$, conditional Jensen's inequality
yields $q_B^2\leq\E[V_1^2\mid B]$, and hence
\begin{align}
\left(\frac1N\sum_{i=1}^N V_i+q_B\right)^2
\leq2\left(\frac1N\sum_{i=1}^N V_i^2+\E[V_1^2\mid B]\right).
\end{align}
The right-hand side is uniformly integrable as $N$ varies. Therefore, \eqref{chaos-path-domination} implies that 
\begin{align}
\{\W_1^{\mathsf X}(Q_B^N,Q_B)^2\}_{N\geq1}
\end{align}
is uniformly integrable. Combining this with the
almost-sure convergence in \eqref{chaos-path-lln}
gives
\begin{align}
\E\W_1^{\mathsf X}(Q_B^N,Q_B)^2
\longrightarrow 0.
\end{align}
Finally, \eqref{chaos-evaluation-bound} yields
\eqref{chaos-sampling-error}.

\smallskip
\noindent\textit{Step 3. Stability of the particles}
For $i=1,\ldots,N$, write
\begin{align}
\Delta Y^i\coloneqq Y^{i,N}-\overline Y^i,\hskip 0.2cm
\Delta Z^i\coloneqq Z^{i,N}-\overline Z^i,\hskip 0.2cm\Delta\mu^i \coloneqq\mu^{i,N}-\overline\mu^i,
\end{align}

and define
\begin{align}\label{chaos-error-norm}
\begin{aligned}
R_t^2&\coloneqq \frac1N\sum_{i=1}^N
\left(|\Delta Y_t^i|^2+
      \norma{\Delta\mu_t^i}_{\BL}^2\right),\\
\mathcal D_N^2
&\coloneqq \E\Bigg[
\sup_{t\in[0,T]}|\Delta\by_t|_N^2
+\int_0^T|\Delta\bz_t|_N^2\,\de t
+\sup_{t\in[0,T]}
\norma{\Delta\bmu_t}_{N,\BL}^2\Bigg].
\end{aligned}
\end{align}
Recall also the particle-copy error $\mathcal E_N^2$
defined in \eqref{normaE}. The distinction is the order of
the supremum and the average: $\mathcal D_N^2$ takes a
supremum of the averaged squared differences, whereas
$\mathcal E_N^2$ averages the individual suprema. Hence
\begin{align}\label{chaos-two-errors}
\mathcal D_N^2\leq\mathcal E_N^2,
\qquad
\E\sup_{t\in[0,T]}R_t^2
\leq\mathcal D_N^2\leq\mathcal E_N^2.
\end{align}
In general, the converse inequality
$\mathcal E_N^2\leq\mathcal D_N^2$ does not hold. Pairing particles and copies with the same index,
as in Step 2, gives
\begin{align}
\W_1(\nu_t^N,\overline\nu_t^N)^2\leq 2R_t^2.
\end{align}
By the triangle inequality and $(a+b)^2\leq 2a^2+2b^2$, one gets
\begin{align}\label{chaos-key-measure-bound}
\begin{aligned}
\W_1(\nu_t^N,m_t)^2&\leq 2\W_1(\nu_t^N,\overline\nu_t^N)^2
+2e_N(t)^2\\
&\leq 4R_t^2+2e_N(t)^2.
\end{aligned}
\end{align}
Let $\Delta h_t^i$ be the difference between the backward drivers of particle $i$ and limit copy $i$. By \autoref{A1}, the inequality $(a+b+c+d)^2\leq4(a^2+b^2+c^2+d^2)$, and \eqref{chaos-key-measure-bound},
\begin{align}\label{chaos-driver-bound}
\begin{aligned}
|\Delta\boldsymbol h_t|_N^2 &\leq
4L_H^2\Bigl(|\Delta\bz_t|_N^2+R_t^2+\W_1(\nu_t^N,m_t)^2\Bigr)\\
&\leq 4L_H^2\Bigl(|\Delta\bz_t|_N^2+5R_t^2+2e_N(t)^2\Bigr).
\end{aligned}
\end{align}
Since each particle and its limiting copy have the same terminal condition and the same initial measure, the corresponding difference equations have zero terminal and initial data. Applying the two estimates of \autoref{lembase} to each pair and averaging over $i$ gives
\begin{align}\label{chaos-backward-stability}
\begin{aligned}
&\E\left[\sup_{t\in[0,T]}|\Delta\by_t|_N^2+\int_0^T|\Delta\bz_t|_N^2\,\de t\right]\\
&\qquad\leq 5T\,\E\int_0^T|\Delta\boldsymbol h_s|_N^2\,\de s\\
&\qquad\leq 20L_H^2T\left[\E\int_0^T|\Delta\bz_s|_N^2\,\de s
+5T\,\E\sup_{t\in[0,T]}R_t^2+2T\varepsilon_N\right].
\end{aligned}
\end{align}
Notice that for the measure components, by subtracting \eqref{chaos-limit-measure} from the particle
equation, the initial measures cancel. Now by applying the argument of Step 5 of the proof of \autoref{importantprop}, and now using \eqref{chaos-key-measure-bound}, one has
\begin{align}\label{chaos-measure-stability}
\begin{aligned}
&\E\sup_{t\in[0,T]}\norma{\Delta\bmu_t}_{N,\BL}^2\\
&\qquad\leq 3L_F^2e^{2L_FT}T
\left[\E\int_0^T|\Delta\bz_s|_N^2\,\de s+5T\,\E\sup_{t\in[0,T]}R_t^2 +2T\varepsilon_N
\right].
\end{aligned}
\end{align}
Since
\begin{align}
\E\sup_{t\in[0,T]}R_t^2\leq\mathcal D_N^2,
\end{align}
combining \eqref{chaos-backward-stability} and \eqref{chaos-measure-stability}, for $T\leq1$, we get
\begin{align}\label{chaos-absorption}
\mathcal D_N^2 \leq C_{\mathrm{st}}T\,\mathcal D_N^2
+C_{\mathrm{err}}T^2\varepsilon_N,
\end{align}
where one may take
\begin{align}
C_{\mathrm{st}} &\coloneqq 120L_H^2+18L_F^2e^{2L_F},\\
C_{\mathrm{err}} &\coloneqq 40L_H^2+6L_F^2e^{2L_F}.
\end{align}
These constants do not depend on $N$. Choose $T_{\mathrm{ch}}\leq1$ so that
\begin{align}\label{chaos-horizon}
C_{\mathrm{fix}}T_{\mathrm{ch}}<1,
\hskip 0.2cm C_{\mathrm{st}}T_{\mathrm{ch}}\leq\frac12,
\end{align}
where $C_{\mathrm{fix}}$ as defined in \eqref{costCfix}. The first condition ensures existence of both the finite and limiting systems, and the second permits absorption in \eqref{chaos-absorption}. Thus
\begin{align}\label{chaos-stability-conclusion}
\mathcal{D}_N^2\leq 2C_{\mathrm{err}}T^2\varepsilon_N
\longrightarrow0.
\end{align}
We next estimate $\mathcal E_N^2$ directly. Apply the two estimates of \autoref{lembase} to each
particle--copy pair, retain the individual supremum in
the estimate for $\Delta Y^i$, and average over $i$.
Apply the Gronwall estimate from Step 5 of the
proof of \autoref{importantprop} in the same way for
$\Delta\mu^i$. Using \eqref{chaos-driver-bound} and
\eqref{chaos-key-measure-bound} gives
\begin{align*}
\mathcal E_N^2
&\leq
\left(20L_H^2+3L_F^2e^{2L_FT}\right)T
\Bigg[
\E\int_0^T|\Delta\bz_s|_N^2\,\de s\\
&\qquad\qquad
+5T\,\E\sup_{t\in[0,T]}R_t^2
+2T\varepsilon_N
\Bigg].
\end{align*}
By \eqref{chaos-two-errors}, the first two error terms
on the right are bounded by $\mathcal E_N^2$.
For $T\leq1$, the constants in
\eqref{chaos-absorption} therefore give
\begin{align}\label{chaos-individual-absorption}
\mathcal E_N^2
\leq C_{\mathrm{st}}T\,\mathcal E_N^2
+C_{\mathrm{err}}T^2\varepsilon_N.
\end{align}
Since $C_{\mathrm{st}}T\leq 1/2$ by
\eqref{chaos-horizon}, we obtain
\begin{align}\label{chaos-individual-stability}
\mathcal E_N^2
\leq 2C_{\mathrm{err}}T^2\varepsilon_N
\longrightarrow0.
\end{align}
\smallskip
\noindent\textit{Step 4. Empirical convergence and finite-dimensional
chaos.}
Since 
\begin{align}
\W_1(\nu_t^N,\overline\nu_t^N) &\leq \frac1N\sum_{i=1}^N
\left(|\Delta Y_t^i|+\norma{\Delta\mu_t^i}_{\BL}\right),
\end{align}
by Cauchy--Schwarz and \eqref{chaos-stability-conclusion},
\begin{align}
\E\sup_{t\in[0,T]}\W_1(\nu_t^N,\overline\nu_t^N)^2
\leq 2\mathcal D_N^2\longrightarrow0.
\end{align}
Applying the triangle inequality with $\overline\nu^N$ and then
\eqref{chaos-sampling-error} proves 
\begin{align}
\begin{aligned}
\E\sup_{t\in[0,T]}\W_1(\nu_t^N,m_t)^2&\leq4\mathcal D_N^2+2\varepsilon_N\\
&\longrightarrow0.
\end{aligned}
\end{align}
Lastly, let us fix $k\in\N$. The coupled systems are exchangeable with respect to the particle indices. Consequently, the expectation of each individual error is the same, and
\eqref{chaos-particle-convergence} implies
\begin{align}\label{chaos-fixed-k}
\begin{aligned}
&\sum_{i=1}^k \E\Bigg[\sup_{t\in[0,T]}|Y_t^{i,N}-\overline Y_t^i|^2+\int_0^T|Z_t^{i,N}-\overline Z_t^i|^2\,\de t\\
&\hskip 0.2cm+ \sup_{t\in[0,T]}
\norma{\mu_t^{i,N}-\overline\mu_t^i}_{\BL}^2
\Bigg]\\
&\qquad=k\mathcal E_N^2\longrightarrow0,
\hskip 0.2cm N\to\infty.
\end{aligned}
\end{align}
In particular, the first $k$ particle paths $(X^{1,N},\ldots,X^{k,N})$ converge in probability, for the product uniform path metric, to $(\overline X^1,\ldots,\overline X^k)$. Conditionally on $B$, the limit paths $\overline X^1,\ldots,\overline X^k$ are independent and have the common conditional path law $Q_B$. Thus their joint conditional law is
\begin{align}
\mathcal L\bigl(\overline X^1,\ldots,\overline X^k\mid B
\bigr)=Q_B^{\otimes k},
\end{align}
and we are done
\end{proof}
Let us now give the proof of \autoref{thm:quantitative-chaos}.
\begin{proof}[{Proof of \autoref{thm:quantitative-chaos}}]
Let $\theta_1,\ldots,\theta_M$ be points supporting the
law of $\Theta^1$, and write
\begin{align}
p_j\coloneqq \P(\Theta^1=\theta_j),\hskip 0.2cm j=1,\ldots,M.
\end{align}
Points with $p_j=0$ can be removed, then we assume $p_j>0$ for every $j$. By the construction of the limit copies in the proof of \autoref{thm:conditional-chaos}, there is a measurable
map of the mark and the Brownian path that produces the entire limit path. Consequently, for each $j$ there is an $\mathsf{X}$-valued random variable $x_j(B)$, measurable with respect to the entire Brownian path, such that
\begin{align}
\overline{X}^i=x_j(B)\hskip 0.2cm\text{on }\{\Theta^i=\theta_j\},
\end{align}
where $\mathsf X=C([0,T],\Y)$. Let us define the empirical frequencies
\begin{align}
\widehat{p}_j^N\coloneqq \frac1N\sum_{i=1}^N
\mathbbm{1}_{\{\Theta^i=\theta_j\}}.
\end{align}
Conditionally on $B$, the probabilities $p_j$ do not change, since the marks are independent of $B$. Thus the conditional path law and the empirical path measure introduced in Part 2 of the previous proof of \autoref{thm:conditional-chaos} have the representations
\begin{align}
Q_B&=\sum_{j=1}^M p_j\delta_{x_j(B)},\\
Q_B^N&=\sum_{j=1}^M\widehat p_j^N\delta_{x_j(B)}.
\end{align}
These representations remain valid if some of the paths $x_j(B)$ coincide. Fix a realization of $B$, and put
\begin{align}
a_j(B) \coloneqq d_{T,\Y}(x_j(B),x_0).
\end{align}
The signed measure $Q_B^N-Q_B$ has total mass zero. Hence, in the Kantorovich--Rubinstein dual representation, a test function $f$ may be replaced by $f-f(x_0)$. For every $1$-Lipschitz $f$,
\begin{align}
\left|\int_{\mathsf X} f\,\de(Q_B^N-Q_B)\right|
&=\left|\sum_{j=1}^M(\widehat p_j^N-p_j)\bigl(f(x_j(B))-f(x_0)\bigr)\right|\\
&\leq\sum_{j=1}^M|\widehat p_j^N-p_j|\,a_j(B).
\end{align}
Taking the supremum over $f$ and applying the weighted Cauchy--Schwarz inequality one gets
\begin{align}\label{chaos-frequency-bound}
\begin{aligned}
\W_1^{\mathsf X}(Q_B^N,Q_B)^2
&\leq \left(\sum_{j=1}^M|\widehat p_j^N-p_j|\,a_j(B)
\right)^2\\
&\leq\left(\sum_{j=1}^M\frac{(\widehat p_j^N-p_j)^2}{p_j}
\right)\left(\sum_{j=1}^M p_j a_j(B)^2\right).
\end{aligned}
\end{align}
Notice that the variables $\mathbbm{1}_{\{\Theta^i=\theta_j\}}$ are i.i.d. Bernoulli variables with parameter $p_j$ and are independent of $B$. Therefore
\begin{align}
\E\left[(\widehat p_j^N-p_j)^2\mid B
\right]=\frac{p_j(1-p_j)}{N}.
\end{align}
Taking conditional expectation in \eqref{chaos-frequency-bound}, we obtain
\begin{align}\label{chaos-conditional-sampling-rate}
\begin{aligned}
\E\left[\W_1^{\mathsf X}(Q_B^N,Q_B)^2\mid B\right]
&\leq\frac1N\left(\sum_{j=1}^M(1-p_j)\right)
\left(\sum_{j=1}^M p_j a_j(B)^2\right)\\
&\leq\frac{M}{N}\int_{\mathsf X}d_{T,\Y}(x,x_0)^2\,Q_B(\de x).
\end{aligned}
\end{align}
Integrating with respect to $B$ and using the definition of $Q_B$ yields
\begin{align}
\E\W_1^{\mathsf X}(Q_B^N,Q_B)^2
&\leq\frac{M}{N}\E\,d_{T,\Y}(\overline X^1,x_0)^2\\
&=\frac{M K_T}{N}.
\end{align}
By the evaluation bound \eqref{chaos-evaluation-bound},
\begin{align}
\varepsilon_N=\E\sup_{t\in[0,T]}\W_1(\overline\nu_t^N,m_t)^2
\leq\E\W_1^{\mathsf X}(Q_B^N,Q_B)^2
\leq\frac{M K_T}{N}.
\end{align}
This proves \eqref{chaos-quantitative-sampling}. By \eqref{chaos-individual-stability} and
\eqref{chaos-quantitative-sampling},
\begin{align*}
\mathcal{E}_N^2
\leq 2C_{\mathrm{err}}T^2\varepsilon_N
\leq \frac{2C_{\mathrm{err}}T^2 M K_T}{N},
\end{align*}
which proves \eqref{chaos-quantitative-particles}. Likewise, the empirical-measure comparison from Part 4 of \autoref{thm:conditional-chaos} one gets
\begin{align}
\E\sup_{t\in[0,T]}\W_1(\nu_t^N,m_t)^2
&\leq 4\mathcal E_N^2+2\varepsilon_N\\
&\leq\frac{(8C_{\mathrm{err}}T^2+2)M K_T}{N},
\end{align}
and thus \eqref{chaos-quantitative-empirical} is proved. Lastly, the coupled particle--copy systems are exchangeable. As observed in \eqref{chaos-fixed-k}, the sum of the errors of the first $k$ indices equals $k\mathcal E_N^2$. Applying \eqref{chaos-quantitative-particles} one has \eqref{chaos-quantitative-fixed-k}, and we are done.
\end{proof}

\section{Application: a population with adaptive strategies}
\label{sec:application}
Consider a population of agents who allocate their activity among
a finite set of strategies. Each agent has a real-valued state
representing the value, or assessment, of an outcome at time $t$.
The outcome is specified at the terminal time $T$, so its earlier
values are determined by a backward equation. At the same time,
the agent adjusts its strategy weights over time in response to
its own state and to information about the population. The
strategy distribution therefore evolves forward, while the state
is determined backward. The two components interact: an agent's
choice of strategies affects its state, and its state affects the
subsequent adjustment of those strategies. Our aim is not to derive an optimal strategy, but to ask whether this coupled description remains meaningful for a large population. In particular, we want to know whether a finite system can be approximated by a representative agent, what role shared uncertainty plays in that approximation, and how large the approximation error is when the population has finitely many
types. The results of the previous sections answer these questions once the model is shown to satisfy their assumptions. Let $d=1$ and let $U=\{1,\ldots,m\}$, endowed with the discrete
metric. An agent's strategy distribution $\lambda\in\cP(U)$ is identified with the vector of weights $(\lambda(1),\ldots,\lambda(m))$. For $\Psi\in\cP_1(\Y)$, define the population observables
\begin{align*}
\overline x(\Psi)&\coloneqq \int_{\Y}x'\,\Psi(\de x',\de\lambda'),\\
s(\Psi)&\coloneqq \int_{\Y}\tanh(x')\,\Psi(\de x',\de\lambda').
\end{align*}
The mean $\overline x(\Psi)$ lets the population's aggregate state enter the backward equation directly. The signal $s(\Psi)$ enters the strategy payoffs instead. Using $\tanh$ makes this signal bounded even when individual states are unbounded, while
preserving its dependence on whether states are positive or
negative. Bounded payoff signals will also make it possible to
verify that the strategy weights remain nonnegative. Fix $a,b,c,r_1,\ldots,r_m\in\R$ and $\alpha_u,\beta_u\in\R$ for $u\in U$. Define the backward driver by
\begin{align}
H(z,(x,\lambda),\Psi)\coloneqq az+bx+c\overline x(\Psi)  +\sum_{u=1}^m r_u\lambda(u).
\label{app-H}
\end{align}
It combines the agent's own state, the population mean, and the
weights assigned to the available strategies. For each strategy
$u$, define the payoff signal
\begin{align}
q_u(x,\Psi)\coloneqq \alpha_u\tanh(x)+\beta_u s(\Psi).
\label{app-payoff}
\end{align}
Thus the payoff signal depends both on the agent's state and on
an aggregate population signal. The strategy weights evolve
according to the replicator rule
\begin{align}
F(z,(x,\lambda),\Psi)(u)\coloneqq \lambda(u)\left(q_u(x,\Psi)-\sum_{v=1}^m\lambda(v)q_v(x,\Psi)\right),\qquad u\in U.
\label{app-F}
\end{align}
A strategy gains weight when its payoff signal exceeds the payoff averaged over the agent's current strategy distribution.
The rule preserves total mass and allows the agent to distribute
its activity among several strategies. Together, \eqref{app-H} and \eqref{app-F} give a coupled model: strategy weights affect the backward state, while that state affects the subsequent adjustment of the weights. Let $(\Theta^i)_{i\geq1}$ be i.i.d. marks, independent of $B$, and choose measurable maps $g:E\to\R$ and $J:E\to\cP(U)$ such that $\E|g(\Theta^1)|^2<+\infty$.
For $\sigma\in\R$, take
\begin{align}
\xi^i=g(\Theta^i)+\sigma B_T, \hskip 0.2cm \overline{\mu}_0^i=J(\Theta^i).
\label{app-data}
\end{align}
The marks allow agents to differ in both their terminal outcomes
and their initial strategy distributions. The term $\sigma B_T$
represents a terminal shock shared by all agents. Consequently,
the large-population limit need not be deterministic: its natural description is the law of a representative agent conditional on the common Brownian motion.
\begin{cor}[Conditional population limit for adaptive strategies]
\label{cor:application}
Consider the $N$-agent system \eqref{eq1} with coefficients
\eqref{app-H} and \eqref{app-F} and data \eqref{app-data}.
There exists $T_{\mathrm{app}}>0$, independent of $N$, such
that for every $0<T\leq T_{\mathrm{app}}$:
\begin{enumerate}
\item[(i)] the $N$-agent system has a unique solution;
\item[(ii)] its empirical measure converges to the conditional
law $m_t$ of a representative agent given the common Brownian
motion, in the sense that
\begin{align}
\E\sup_{t\in[0,T]}\W_1(\nu_t^N,m_t)^2
\longrightarrow0;
\label{app-empirical-limit}
\end{align}
\item[(iii)] if the law of $\Theta^1$ is supported on at most
$M$ points, then, with $K_T$ and $C_{\mathrm{err}}$ as in
\autoref{thm:quantitative-chaos},
\begin{align}
\mathcal E_N^2&\leq \frac{2C_{\mathrm{err}}T^2 M K_T}{N},
\label{app-particle-rate}\\
\E\sup_{t\in[0,T]}\W_1(\nu_t^N,m_t)^2&\leq \frac{(8C_{\mathrm{err}}T^2+2)M K_T}{N}.
\label{app-empirical-rate}
\end{align}
\end{enumerate}
\end{cor}

\begin{proof}
We verify the assumptions of \autoref{thm:conditional-chaos} and
\autoref{thm:quantitative-chaos}. The maps $x\mapsto x$ and $x\mapsto\tanh(x)$ are $1$-Lipschitz. By the Kantorovich--Rubinstein inequality, both $\overline{x}$ and $s$ are therefore Lipschitz with respect to $\W_1$. Moreover, $|s(\Psi)|\leq1$.
Since $U$ is finite, the term $\sum_u r_u\lambda(u)$ is Lipschitz in $\lambda$ for the $\BL$ norm. This proves \autoref{A1} for \eqref{app-H}. Let us set

\begin{align}
Q\coloneqq\max_{u\in U}(|\alpha_u|+|\beta_u|).    
\end{align}
Then $|q_u(x,\Psi)|\leq Q$. Each $q_u$ is Lipschitz in $(x,\Psi)$, and on the finite space $U$ the coordinate maps $\lambda\mapsto\lambda(u)$ are Lipschitz for the $\BL$ norm. Expanding the products in \eqref{app-F} consequently shows that $F$ is globally Lipschitz in $(x,\lambda,\Psi)$. It does not depend on $z$, so \autoref{A2} follows. Let
\begin{align}
\overline q(x,\lambda,\Psi)
\coloneqq \sum_{v=1}^m \lambda(v)q_v(x,\Psi).
\end{align}
Since $\sum_{u=1}^m\lambda(u)=1$, 
\begin{align}
\sum_{u=1}^m F(z,(x,\lambda),\Psi)(u)
&=\sum_{u=1}^m\lambda(u)
\bigl(q_u(x,\Psi)-\overline q(x,\lambda,\Psi)\bigr)\notag\\
&=\overline q(x,\lambda,\Psi)
-\overline q(x,\lambda,\Psi)\sum_{u=1}^m\lambda(u)\notag\\
&=0.
\end{align}
Thus \autoref{A3} holds. Furthermore, for every $u\in U$,
\begin{align*}
F(z,(x,\lambda),\Psi)(u)+2Q\lambda(u)
&=\lambda(u)\left(q_u(x,\Psi)-\sum_{v=1}^m\lambda(v)q_v(x,\Psi)+2Q
\right)\\
&\geq0.
\end{align*}
Thus \autoref{A4} holds with $\delta=2Q$. The claims follow
from \autoref{cor:random-initial} (applied with the enlarged
filtration \eqref{chaos-filtration} and the random initial data
\eqref{chaos-data}), \autoref{thm:conditional-chaos}, and
\autoref{thm:quantitative-chaos}.
\end{proof}
The corollary identifies the population limit in this model. Even for a large population, the limiting law $m_t$ remains
random because it depends on the common Brownian motion. Conditional on that noise, individual agents become independent
in the limit. 

When the law of the agent marks has finite support,
\eqref{app-particle-rate} and \eqref{app-empirical-rate}
quantify the finite-population approximation: the squared particle and empirical-measure errors are of order $N^{-1}$.
\subsection{Numerical illustration}
\label{subsec:app-numerics}
We illustrate Corollary~\ref{cor:application} with two strategies
and three agent types. The type probabilities, terminal values,
and initial strategy distributions are, respectively,
\begin{align*}
p&=(0.25,\,0.45,\,0.30),\hskip 0.2cm g=(-1,\,0.2,\,1.2),\\
\lambda_{0,1}&=(0.80,\,0.20),\hskip0.2cm \lambda_{0,2}=(0.50,\,0.50),\hskip 0.2cm \lambda_{0,3}=(0.15,\,0.85).
\end{align*}
We take
\begin{align*}
&T=0.2,\qquad a=0.3,\hskip 0.2cm b=0.4,\hskip 0.2cm  c=0.5,
\hskip 0.2cm  \sigma=0.7,\hskip 0.2cm r_1=r_2=0,\\
&(\alpha_1,\alpha_2)=(0.8,-0.6),\hskip 0.2cm 
(\beta_1,\beta_2)=(-0.3,0.5).
\end{align*}
The choice $T=0.2$ is made for numerical illustration.
It is not covered by the explicit sufficient time-horizon
condition used in the proof of \autoref{thm:conditional-chaos} for these parameters. Consequently, the experiment illustrates the observed finite-population behaviour, rather than providing a
numerical instance within the time horizon guaranteed by
our theorems. 

We set $r_1=r_2=0$ to obtain an explicit solution for the
backward component. This choice removes the feedback from
strategy weights to the backward state in the numerical example:
the states affect the replicator dynamics, but the strategy
weights do not affect the states. Thus the experiment illustrates
the finite-population errors in a one-way-coupled special case
of Corollary~\ref{cor:application}, rather than numerically
testing the fully coupled model.

Write $\overline g=\sum_{j=1}^3p_jg_j$. In a finite sample, let $\widehat p_j^N$ be the proportion of agents of type $j$ and put $\overline g_N=\sum_{j=1}^3\widehat p_j^Ng_j$. For an agent of type $j$, the backward solutions in the conditional limit and in the finite system are
\begin{align}
\overline Y_j(t)
&=e^{(b+c)(T-t)}
   \bigl(\overline g+\sigma B_t+a\sigma(T-t)\bigr)
  +e^{b(T-t)}(g_j-\overline g), \label{app-num-limit}\\
Y_j^N(t)&=e^{(b+c)(T-t)}\bigl(\overline g_N+\sigma B_t+a\sigma(T-t)\bigr)+e^{b(T-t)}(g_j-\overline g_N). \label{app-num-finite}
\end{align}
Indeed, in the finite system, the empirical mean $M_t^N=\sum_{j=1}^3\widehat p_j^N Y_j^N(t)$ satisfies a backward equation with coefficient $b+c$ and terminal value $\overline g_N+\sigma B_T$. Subtracting this mean equation from the equation for type $j$ leaves a deterministic backward equation for $Y_j^N-M^N$, with coefficient $b$. The same argument applies to the conditional limit, with $p_j$ in place of $\widehat p_j^N$. In both systems,
\begin{align*}
Z_j(t)=\sigma e^{(b+c)(T-t)}.
\end{align*}
Thus the particle--limit error in this example has no $Z$
contribution. In particular,
\begin{align*}
Y_j^N(t)-\overline Y_j(t)=\bigl(e^{(b+c)(T-t)}-e^{b(T-t)}\bigr)(\overline g_N-\overline g).
\end{align*}
This identity makes explicit how fluctuations in the sampled
type proportions affect the finite-population states. The strategy distributions are obtained by applying the
replicator equation \eqref{app-F} to the three types. We use $100$ equal time steps, with $\Delta t=T/100$. In each repetition, the types of the $N$ agents are sampled independently with probabilities $p$, and a single common Brownian path is generated from independent increments $\Delta B_n=\sqrt{\Delta t}\,G_n$, where the $G_n$ are standard Gaussian. The finite system and the conditional limit use the same path within that repetition. Their backward states are evaluated at the grid times using \eqref{app-num-limit} and \eqref{app-num-finite}. For the strategy distributions, we use explicit Euler updates. At each grid time, the payoff signals are computed from the current strategy-independent states. The population
signal is weighted by $\widehat p^N$ in the finite system and by $p$ in the conditional limit; each strategy vector is then advanced using $\Delta t$ times the right-hand side of \eqref{app-F}. The MATLAB script uses the random seed \texttt{rng(7)}. For clarity, we specify the metric used in the numerical Wasserstein calculation. If $U=\{1,2\}$ has the discrete metric and $\lambda,\eta\in\cP(U)$, the definition
\eqref{normLip} gives
\begin{align}
\|\lambda-\eta\|_{\BL}=\frac23|\lambda(1)-\eta(1)|.
\label{app-two-strategy-BL}
\end{align}
To see this, put $v=\phi(1)-\phi(2)$. For fixed $v$, the
smallest possible value of $\|\phi\|_\infty$ is $|v|/2$.
Hence the constraint $\|\phi\|_\infty+\mathrm{Lip}(\phi)\leq1$ permits $|v|\leq2/3$, and this bound is attained. The transport
cost between two state--strategy pairs is consequently
\begin{align}
d_{\Y}\bigl((x,\lambda),(y,\eta)\bigr)=|x-y|+\frac23|\lambda(1)-\eta(1)|.
\label{app-numerical-cost}
\end{align}
At a grid time $t_n$, the empirical measure and conditional
limit law can be written as
\begin{align*}
\nu_{t_n}^N&=\sum_{j=1}^3\widehat p_j^N
\delta_{(Y_j^N(t_n),\lambda_j^N(t_n))},\\
m_{t_n}&=\sum_{j=1}^3p_j\delta_{(\overline{Y}_j(t_n),\overline\lambda_j(t_n))}.
\end{align*}
We compute their Wasserstein distance by solving the
three-by-three transport problem
\begin{align}
\W_1(\nu_{t_n}^N,m_{t_n})
=\min_{\substack{\pi_{jk}\geq0\\
                  \sum_k\pi_{jk}=\widehat p_j^N\\
                  \sum_j\pi_{jk}=p_k}}
\sum_{j,k=1}^3\pi_{jk}\,
d_{\Y}\bigl((Y_j^N(t_n),\lambda_j^N(t_n)),
(\overline Y_k(t_n),\overline\lambda_k(t_n))\bigr).
\label{app-numerical-transport}
\end{align}
In the implementation, one redundant marginal constraint is removed, leaving five independent equality constraints on the nine transport variables. All nonsingular five-column bases of the resulting constraint matrix are precomputed. At each grid time, the code evaluates their basic solutions, retains those with nonnegative flows up to a numerical tolerance of $10^{-10}$, and takes the smallest transport cost. For each population size $N$, the experiment is repeated
$300$ times with independently sampled types and Brownian paths. If $t_0,\ldots,t_{100}$ is the simulation grid, the
blue curve is the Monte Carlo estimate of
\begin{align}
\sum_{j=1}^3\widehat p_j^N\left[\max_n|Y_j^N(t_n)-\overline Y_j(t_n)|^2+\max_n\|\lambda_j^N(t_n)-\overline\lambda_j(t_n)\|_{\BL}^2\right],
\label{app-numerical-particle-error}
\end{align}
averaged over repetitions. The orange curve is the Monte Carlo estimate of
\begin{align}
\max_n\W_1(\nu_{t_n}^N,m_{t_n})^2,
\label{app-numerical-wasserstein-error}
\end{align}
also averaged over repetitions. These are time-grid
approximations of the two errors in \eqref{app-particle-rate} and \eqref{app-empirical-rate}.

\begin{figure}[ht!]
    \centering
    \includegraphics[width=0.6\textwidth]{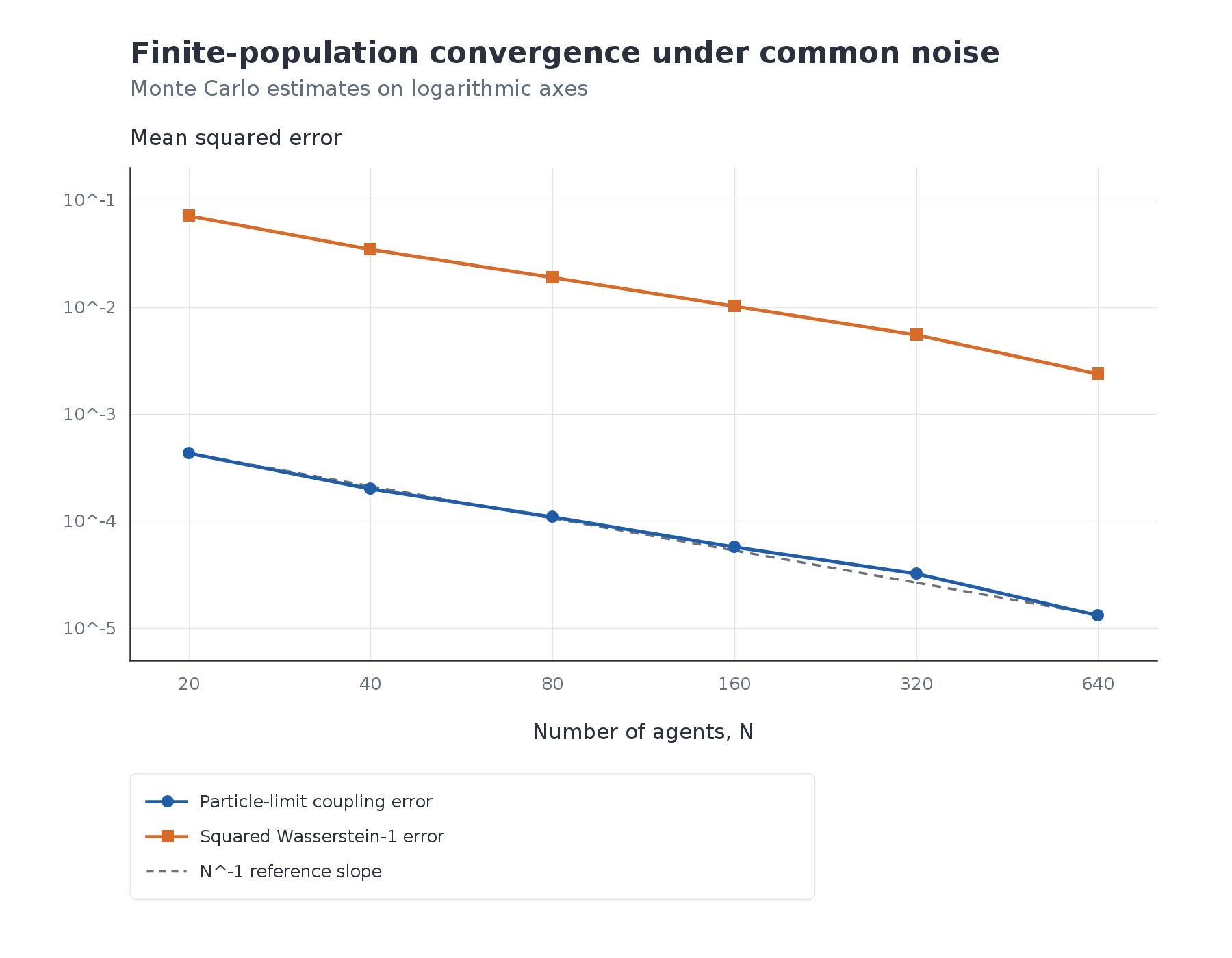}
    \caption{Monte Carlo estimates of the particle--limit error
    \eqref{app-numerical-particle-error} and the squared
    empirical-measure Wasserstein error
    \eqref{app-numerical-wasserstein-error}, plotted against
    population size $N$. Both systems share the same common
    Brownian path within each repetition. The dashed line is an
    $N^{-1}$ reference slope, not a theoretical upper bound.}
    \label{fig:adaptive-strategies}
\end{figure}
Both errors decrease with the population size in
Figure~\ref{fig:adaptive-strategies}, with an observed
squared-error decay close to the $N^{-1}$ reference slope.
The Wasserstein error also reflects fluctuations of the
empirical type frequencies $\widehat p^N$ around $p$,
whereas the particle error compares each agent only with
its limit copy of the same type. This helps explain the
different magnitudes of the two curves. Since $T=0.2$ is
outside the time horizon guaranteed by the sufficient
condition in our proof, the figure should not be read as
a verification of the quantitative bounds in
\autoref{cor:application} for these parameters.

\section{Conclusion}
We have studied an interacting system in which each agent has a backward stochastic state and a probability-valued strategy distribution evolving forward in time. Both components depend on the empirical state--strategy law, and all agents are driven by a common Brownian motion. This combination requires the backward
equation and the probability-valued evolution to be treated as a coupled system: the former is prescribed at the terminal time, whereas the latter is prescribed at the initial time and must remain a probability measure throughout its evolution.

Under Lipschitz, mass-conservation, and positivity assumptions, we proved well-posedness on a sufficiently short time interval, with a time horizon independent of the number of agents. The fixed-point construction combines a backward equation obtained through Brownian martingale representation with a positivity-preserving construction for the forward measure-valued equation. The resulting estimates are uniform in the population size and also apply when independently sampled agent marks are included in the initial filtration.

For the large-population limit, we identified the relevant interaction law as the law of a representative state--strategy pair conditional on the common Brownian motion. We proved convergence of the empirical law to this conditional law, uniformly in time in expected squared Wasserstein distance. The limiting agent paths are independent conditionally on the common Brownian path, although they generally remain dependent without conditioning. Thus the common noise persists in the
population limit rather than being averaged out.

When the distribution of agent marks has finite support, the conditional law of the limiting paths is supported on finitely many type-dependent paths. In this setting, fluctuations of the sampled type frequencies yield an $N^{-1}$ bound for the
expected squared sampling error. Stability estimates transfer this bound to the average squared particle--copy error and to the expected squared Wasserstein distance between the empirical and conditional limiting laws. These quantitative conclusions concern squared errors and rely on the finite-type assumption; they do not assert the same rate for general mark
distributions.

Two limitations delimit the present theory. First, the
contraction and stability arguments require a sufficiently short time horizon; extending them to longer intervals would require additional structural conditions or a different method.
Second, quantitative bounds for more general mark distributions would require control of the conditional empirical sampling error beyond the finite-support argument used here. Investigating these extensions are natural next steps.

\section*{Acknowledgements}
The author thanks the Fields Institute for Research in Mathematical Sciences for its hospitality during the preparation of this work.
\bibliographystyle{siam}
\bibliography{bibliography}

\end{document}